\documentclass[12pt,a4paper]{amsart}

\usepackage{amssymb}
\let\g=\mathfrak

\let\dag=\perp
\newcommand{\0}{{\boldsymbol{0}}}
\newcommand{\1}{{\boldsymbol{1}}}
\newcommand{\ii}{{\boldsymbol{i}}}
\newcommand{\bb}{{\boldsymbol{b}}}

\newcommand{\bw}{{\boldsymbol{w}}}
\newcommand{\bB}{{\boldsymbol{B}}}

\newcommand{\Hom}{\operatorname{Hom}}
\newcommand{\Spf}{\operatorname{Spf}}
\newcommand{\Pfaff}{\operatorname{Pfaff}}
\newcommand{\Ber}{\operatorname{Ber}}
\newcommand{\tr}{\operatorname{tr}}
\newcommand{\str}{\operatorname{str}}
\newcommand{\singsupp}{\operatorname{singsupp}}

\def\AC{{\mathcal A}}

\def\CC{{\mathcal C}}
\def\DC{{\mathcal D}}

\def\FC{{\mathcal F}}

\def\HC{{\mathcal H}}
\def\IC{{\mathcal I}}
\def\JC{{\mathcal J}}

\def\PC{{\mathcal P}}
\def\QC{{\mathcal Q}}

\def\SC{{\mathcal S}}

\def\UC{{\mathcal U}}
\def\VC{{\mathcal V}}
\def\WC{{\mathcal W}}
\def\XC{{\mathcal X}}

\let\a=\alpha
\let\be=\beta

\let\C=\Gamma

\let\e=\epsilon

\let\l=\lambda
\let\L=\Lambda
\let\w=\omega
\let\W=\Omega
\let\r=\rho
\let\s=\sigma
\let\d=\delta

\let\de=\partial

\newtheorem{theo}{Theorem}[section]
\newtheorem{defi}{Definition}[section]
\newtheorem{coro}{Corollary}[section]
\newtheorem{prop}{Proposition}[section]
\newtheorem{lemme}{Lemma}[section]

\let\what=\widehat
\def\som{\mathop{\sum}\limits}

\def\union{\mathop{\cup}\limits}
\def\tens{\mathop{\otimes}\limits}
\def\sdir{\mathop{\oplus}\limits}

\def\sup{\mathop{Sup}\limits}

\let\oo=\infty

\def\hfl#1#2{\smash{\mathop{\hbox to
12mm{\rightarrowfill}}\limits^{\scriptstyle#1}_{\scriptstyle#2}}}

\def\ihfl#1#2{\smash{\mathop{\hbox to
12mm{\hookrightarrowfill}}\limits^{\scriptstyle#1}_{\scriptstyle#2}}}

\def\lihfl#1#2{\smash{\mathop{\hbox to
12mm{\hookleftarrowfill}}\limits^{\scriptstyle#1}_{\scriptstyle#2}}}

\def\lhfl#1#2{\smash{\mathop{\hbox to
12mm{\leqslantftarrowfill}}\limits^{\scriptstyle#1}_{\scriptstyle#2}}}

\title{superPfaffian}
\author{ P. Lavaud}

\begin{document}

\begin{abstract} Let  $V=V_\0\oplus V_\1$ be a  real finite dimensional
supervector space  provided with a non-degenerate antisymmetric even bilinear
form $B$. Let  $\g{spo}(V)$  be the Lie superalgebra of endomorphisms of $V$
which preserve $B$. We consider  $\g{spo}(V)$ as a supermanifold. We show that a
choice of an orientation of $V_\1$ and of a square root $\ii$ of $-1$ determines
a very interesting generalized function on the supermanifold $\g{spo}(V)$, the
\emph{superPfaffian}.

When $V=V_\1$, $\g{spo}(V)$ is the orthogonal Lie algebra $\g{so}(V_\1)$ and the
superPfaffian is the usual  Pfaffian, a square root of the determinant.

When $V=V_\0$, $\g{spo}(V)$ is the symplectic Lie algebra $\g{sp}(V_\0)$ and the
superPfaffian  is a constant multiple of  the Fourier transform of one the two
minimal nilpotent orbits in the dual of the Lie algebra $\g{sp}(V_\0)$, and is
an analytic  square root of the inverse of the determinant in the open subset of
invertible elements of $\g{spo}(V)$.

Our opinion is that the superPfaffians (there are four of them, corresponding to
the two   orientations on $V_\1$, and to the two square roots of $-1$) are 
fundamental objects. At least, they occur in the study of equivariant cohomology
of supermanifolds (\cite{Lav98,Lav04}), and in the study of the metaplectic
representation of the metaplectic group 
with Lie algebra $\g{spo}(V)$. In this article, we present the definition and
some basic properties of the superPfaffians.
\end{abstract}

\maketitle

\section*{Introduction} Let $V$ be an oriented  finite dimensional real   
vector space provided with a non degenerate symmetric bilinear form $B$.
 On $\g{so}(V)$, the Pfaffian is a polynomial square root of the determinant,
and it is well known that the Pfaffian of $X\in\g{so}(V)$ can be defined by a
suitable Berezin integral (the notation $\int_V d_V$ for the integral is
defined  in section \ref{integration})
\begin{eqnarray}
\label{eq:1}
 \int_V d_V(v) \exp (-\frac12 B(v,Xv) )
\end{eqnarray}
over the vector space $V$ seen as an odd space (cf. for example
\cite{BGV92} and section \ref{analytic} below).

\medskip

This definition still  have  a formal meaning for a real supervector space
$V=V_\0\oplus V_\1$ provided with a non-degenerate antisymmetric even bilinear
form $B$ and an orientation of $V_\1$: $V_\0$ is a symplectic vector space,
$V_\1$ is an oriented vector space provided with a non degenerate symmetric
bilinear form, and this structure provides us with a well defined
\emph{Liouville integral} $d_V$ on the supermanifold $V$(it is a specific
normalization of the Berezin integral on $V$). However, the integral
(\ref{eq:1}) is convergent only when $X$ is in  an open subset of
$\g{spo}(V)$ (definition is given in section \ref{analytic} below),
where it is (as expected, and well known), a square root of the inverse
Berezinian (Definition is given in section  \ref{sec:Ber} below). We
call this function the superPfaffian. For instance, if $V=V_0$, then for
$X\in\g{sp}(V)$, $B(v,Xv)$ is a quadratic form on $V$, and the integral
(\ref{eq:1}) is convergent when it is positive definite. In this case, $\det(X)$
is strictly positive, and the superPfaffian is the positive square root of
$1/\det(X)$.

\medskip Since the inverse Berezinian is not polynomial (it is only a rational
function when $V_0$ is not $0$), there is no natural extension of this
superPfaffian to a function on the supermanifold $\g{spo}(V)$. The purpose of
this article is to show that  there is a natural extension of this superPfaffian
as a
\emph{generalized function} on the supermanifold $\g{spo}(V)$. Notice that the
superPfaffian is $0$ when $\dim(V_\1)$ is odd. Let $m=\dim(V_0)$ (which is even)
and $n=\dim(V_\1)$ (which we assume now to be even). Let $\ii\in\mathbb{C}$ be a
square root of $-1$. We define the superPfaffian   by the formula:
\begin{eqnarray}
\label{eq:2}
\ii^{(m-n)/2} \int_V d_V(v) \exp (-\frac\ii 2 B(v,Xv) ).
\end{eqnarray}
Of course we   prove that  (\ref{eq:2}) has a well defined meaning as a
generalized function of $X$ on the supermanifold $\g{spo}(V)$.

\bigskip

The superPfaffian has very nice properties. As in \cite{MQ86} for $V=V_{\1}$, some of them 
follow from the transformations properties of $d_{V}$ under linear or affine change of variables.

a) It is an analytic function in the open set where  the inverse Berezinian is
defined, and, in this open set,  it is a square root of the inverse Berezinian.

b) In the open set  where the   integral (\ref{eq:1}) is convergent, it is equal
to the function given by (\ref{eq:1}).

c) It is a boundary value of an holomorphic function defined in a specific cone
of $\g{spo}(V\otimes  \mathbb C)$

d) It is harmonic. More precisely, let $SpO(V)$ be the supergroup of
endomorphisms of $V$ which preserve $B$.  The superPfaffian is annihilated by
the homogeneous constant coefficient differential operators on $\g{spo}(V)$
which are of degree $>0$ and $SpO(V)$-invariant.

\medskip

We prove that properties a) and c)  determine the superPfaffian up to sign.

\medskip For some values of $m$ and $n$, we are able to prove that properties a)
and d) determine  the superPfaffian up to sign and complex conjugation. However,
we do not know if this is true in general.

\bigskip One main motivation for this study comes from differential
supergeometry. The (equivariant) Euler form of a (equivariant) real Euclidean
oriented   fiber bundle is equal  to the Pfaffian of the (equivariant) curvature
of an (equivariant) connection on this bundle. With complications partly due to
the fact that the superPfaffian is not a function but only a generalized
function, this is still true in supergeometry (cf
\cite{Lav98,Lav02}). Thus the superPfaffian plays an important role in the
equivariant cohomology of supermanifolds, in particular in relation with  the
localization formula. In fact, formulas (\ref{eq:1}) and (\ref{eq:2})  may be
considered as particular typical cases of the localization formula.

\medskip A second  motivation is the close relationship  between the
superPfaffian and the distribution character of the metaplectic representation
of the simply  connected Lie supergroup with Lie superalgebra $\g{spo}(V)$. This
will be studied in another paper.

\medskip 

I wish to thank Michel Duflo who introduced me to supermathematics and spent much
of his time 
to suggest to me many deep  improvements to this paper.

\newpage

\tableofcontents

\section{Prerequisites}
\subsection{Notations} In this article, unless otherwise specified, all
supervector spaces and superalgebras will be real.  If $V$ is a supervector
space, we denote by $V_\0$ its even part and by $V_\1$ its odd part. If $v$ is a
non zero homogeneous element of $V$, we denote by $p(v)\in \mathbb Z/2\mathbb Z$
its parity.   We put $\dim(V)=(\dim(V_\0),\dim(V_\1))$. We denote by $V^*$ the
dual supervector space $\Hom(V,\mathbb R)$. If $V$ and $W$ are supervector
spaces, $V\otimes W $ and  $W\otimes V$ are supervector spaces, and they are
identified using the rule of signs (for non zero homogenous $v\in V$ and $w\in
W$ we identify $v\tens w$ and 
$(-1)^{p(v)p(w)}w\tens v$). We denote by $S(V)$ the symmetric algebra of $V$.
Recall that it is equal to $S(V_\0)\otimes
\Lambda(V_\1)$, where $S(V_\0)$ and $\Lambda(V_\1)$ are the classical symmetric
and exterior algebras of the corresponding ungraded vector spaces. We  
use the notation $\Lambda(U)$ only for ungraded vector spaces $U$. So, if $V$ is
a supervector space, $\Lambda(V)$ is the exterior algebra of the underlying
vector space.

Let $(m,n)\in \mathbb N\times \mathbb N$. We denote by $\mathbb R^{(m,n)}$ the
supervector space of dimension $(m,n)$ such that $V_\0=\mathbb R^m$ and
$V_\1=\mathbb R^n$.

We choose  a square root $\ii$ of $-1$.

\subsection{Near superalgebras} We say that a commutative superalgebra $\PC$ is
\emph{near} 
 if it is finite dimensional, local, and  with $\mathbb R$ as
residual field. They are the \emph{alg\`{e}bres proches} of Weil
\cite{Wei53}. For $\alpha \in\PC$, we denote by $\bb(\alpha )$ the canonical
projection of $\alpha $ in $\mathbb R$ ($\bb(\alpha )$ is the \emph{body} of
$\alpha $, and $\alpha -\bb(\alpha )$ ---a nilpotent element of $\PC$--- the
\emph{soul} of $\alpha $, according to the terminology of \cite{DeW84}). Let
$\a\in\PC_\0$ be an even element. Let $\phi\in\CC^{\oo}(\mathbb{R},W)$ be a
smooth function defined in a neighborhood of $\bb(\alpha )$ in $\mathbb R$, with
values in some Fr\'echet supervector space $W$. We freely use the notation:
\begin{equation}
\label{FonctLisse} \phi(\alpha )=\sum_{k=0}^\infty \frac{(\alpha -\bb(\alpha
))^k }{k!}\phi^{(k)}(\bb(\alpha ))\in W\otimes \PC.
\end{equation}
In particular, if $\alpha\in \PC_\0$ is invertible, its absolute value $|\alpha
|\in \PC_\0$ is defined by the formula:
\begin{equation}
|\alpha |={\frac{|\bb(\alpha )|}{ \bb(\alpha )}}\alpha ,
\end{equation}
and if $\bb(\a)>0$, its square root is defined by the finite sum 
\begin{displaymath}
\sqrt\a=\sqrt{\bb(\a)}\Big(1+\frac12(\frac{\a}{\bb(\a)}-1)-\frac
1{2^22!}(\frac{\a}{\bb(\a)}-1)^2+\frac 3 {2^3 3!}(\frac{\a}{\bb(\a)}-1)^3 -\frac
{3.5} {2^4 4!}(\frac{\a}{\bb(\a)}-1)^4+\dots\Big),
\end{displaymath}
 where $\sqrt \l$ is the unique positive square root of $\l>0$. A notation like
$\sqrt{|\a|}$ (for $\a\in\PC_{\0}$ invertible) is not ambiguous, because if
$f\in\CC^\oo(\mathbb{R}^+)$ and $\g g\in\CC^{\oo}(\mathbb R,\mathbb R^{+})$ (where $\mathbb R^+=\{x\in\mathbb R,\,x\ge 0\}$),
 for $\a\in\PC_{\0}$ $f\circ g(\a)= f(g(\a))$.

\subsection{Supermanifolds}\label{superman} By a \emph{supermanifold} we mean a
smooth real supermanifold as in
\cite{Kos77}, \cite{Ber87}, \cite{BL75}. Let $V$ be a finite dimensional
supervector space. We denote also by $V$ the associated supermanifold. In this
paper we mostly use this kind of supermanifolds, and we recall some relevant 
definitions in this particular case.

Let $\UC \subset V_\0$ be an open set. We put $$
\CC_V^\infty(\UC)=\CC^\infty(\UC)\otimes \Lambda(V_\1^*),$$
 where $\CC^\infty(\UC)$ is the usual algebra of smooth real valued functions
defined in $\UC$, and $\Lambda(V_\1^*)$ is the exterior algebra of $V_\1^*$. We
say that $\CC_V^\infty(\UC)$ is the superalgebra of \emph{smooth functions on
$V$ defined in $\UC$}. The supermanifold $V$ is by definition the topological
space $V_\0$ equipped with the sheaf of superalgebras $\CC_V^\infty$. Notice
that if $\UC$ is not empty, there is a canonical inclusion $S(V^*)
\subset\CC^\infty(\UC)$. The corresponding elements are called \emph{polynomial
functions}. One can also define rational functions. Similarly, if $W$ is a
Fr\'{e}chet supervector space (for instance $W=\mathbb C$), we denote by
 $\CC_V^\infty(\UC,W)=\CC^\oo(\UC,W)\tens\L(V_1)^*$ the space of $W$-valued
smooth functions.

\medskip

Let $\PC$ be  a near superalgebra. We put $$
V_\PC=(V\otimes \PC)_\0. $$
It is called the \emph{set of points of $V$ with values in $\PC$}. Extending the
body $\bb :\PC \to \mathbb R$ to a map $V\otimes \PC
\to V$, and restricting it  to $V_\PC$, we obtain a map, also called the body
and denoted by $\bb$, $$
\bb :V_\PC \to V_\0. $$
Let  $\UC \subset V_\0$ be an open set. We denote by $V_\PC(\UC)
\subset V_\PC$ the inverse image of $\UC$ in $V_\PC$. It is known that
$V_\PC(\UC)$ is canonically identified to the set of (even) algebra
homomorphisms   $\CC_V^\infty(\UC) \to \PC$. Let $v\in V_\PC(\UC)$. We will
denote the corresponding character by $\phi
\mapsto \phi(v)$, and say that $\phi(v) \in \PC$ is the value of $\phi\in
\CC_V^\infty(\UC) $ at the point $v$.

\medskip For example, let $ v=v_i p^i \in V_\PC$ (we use Einstein's summation
rule, and considering tensorisation by $\PC$ as an extension of scalars, we
write $v_i p^i$ instead of $v_i \otimes p^i$) where the $(v_i,p^i)\in V\times
\PC$ are a finite number of pair of homogeneous elements with the same parity.
Then $$
\bb(v)=v_i \bb(p^i) , $$
 which is in $V_\0$ since $\bb(p^i)=0$ if $p^i$ is odd. Let $\phi\in V^*$. We
denote by the same letter the corresponding element in $\CC_V^\infty(V_{0})$.
Then $\phi(v)=\phi(v_i)p^i$, and this formula in fact completely determines the
bijection between $V_\PC(\UC)$ and the set of even homomorphisms of algebras 
$\CC_V^\infty(\UC) \to \PC$ (cf.
\cite{Wei53}).

\medskip For   $\phi \in \CC_V^\infty(\UC)$, we denote by $\phi_\PC$ the
corresponding function $v \mapsto \phi(v)$ defined in $V_\PC(\UC)$. Then
$\phi_\PC \in \CC^\infty(V_\PC(\UC),\PC)$. The importance of this construction
is that for $\PC$ large enough (for instance if $\PC$ is an exterior algebra
$\Lambda \mathbb R^N$ with $N\geq
\dim(V_\1)$), the map $\phi \mapsto \phi_\PC$ is injective, which allows more or
less to treat   $\phi$ as an ordinary function.

\medskip We emphasize  the special case $\PC=\mathbb R$. Then $V_\mathbb
R=V_\0$, $V_{\mathbb R}(\UC)=\UC$, and $\phi_\mathbb R$ is the projection of
$\phi\in \CC_V^\infty(\UC)$ to $\CC^\infty(\UC)$ which naturally extends the
projection of $\Lambda(V_\1^*)$ to $\mathbb R$.

\medskip To help the reader, we give two typical examples.

\medskip
 Let  $V=\mathbb R^{(1,0)}$. Then $V_\0=V=\mathbb R$, and $V_\PC= \PC_\0$. Let
$\UC \subset\mathbb R$ be an open subset,
$\phi\in\CC_V^\infty(\UC,W)=\CC^\infty(\UC,W)$, and $\alpha\in V_\PC=
\PC_\0$ such that $\bb(\alpha)\in \UC$. Then $\phi(\alpha)\in \PC$ is defined by
formula (\ref{FonctLisse}).

\medskip
 Let   $V=\mathbb R^{(0,1)}$.  Then $V_\0=\{0\}$, and $V_\PC=
\PC_\1$.  Any $\phi\in\CC_V^\infty( \{0\},W) $ can be written as $\phi=c+\xi d$,
where $c$ and $d$ are elements of $W$, and $\xi$ the standard coordinate (the
identity function) on $\mathbb R$. Then, for $\alpha \in  \PC_\1$,
$\phi(\alpha)\in \PC$ is defined by formula $$
\phi(\alpha)=c+\alpha d .$$

\medskip Let $\AC$ be a commutative superalgebra. We still use the notation
$V_\AC=(V\otimes \AC)_\0$.  Polynomial functions $S(V^*)$ can be evaluated on
$V_\AC$, but,  in general, smooth functions can be evaluated  on $V_\AC$ only if
$\AC$ is a near algebra.

\medskip The particular case $\AC=S(V^*)$ is important, because $V_\AC$ contains
a particular point, the \emph{generic point}, corresponding to the identity in
the identification of $Hom(V,V)_{\0}$ with $(V\otimes V^* )_{\0}\subset V_\AC$.
Let us call $v$ the generic point. Then we have $f(v)=f$ for any polynomial
function $f \in S(V^*)$.

\subsection{Coordinates and integration}\label{integration}

 Let $V$ be a finite dimensional supervector space. By  a  basis $(g_i)_{i\in
I}$ of $V$, we mean an indexed basis consisting of homogeneous elements. The  
\emph{dual basis} $(z^i)_{i\in I}$ of $V^*$ is defined by the usual relation
$z^j(g_i)=\delta^j_i$ (the Dirac symbol). We will also say that the basis
$(g_i)_{i\in I}$ is  the
\emph{predual basis} of the basis $(z^i)_{i\in I}$ (the dual basis, in the
canonical identification of $V$ to the dual 
of $V^*$ is  $((-1)^{p(g_i)}g_i)_{i\in I}$). A basis $(z^i)_{i\in I}$ of $V^*$
will be also called \emph{a system of coordinates on $V$}. The corresponding
vector fields (i.e. derivations of the algebra of smooth functions on $V$) are
denoted by $\frac{\de}{\de z^i}$. They are characterized by the rule $$
\frac{\de}{\de z^j}(z^i)=\delta^i_j. $$
Notice that the generic point $v$ of $V$ is then given by the formula
\begin{eqnarray}
\label{eq:generic} v=g_i z^i \in V_{S(V^*)}.
\end{eqnarray}

\medskip We will mainly use standard coordinates. Let $(m,n)=\dim(V)$. Then they
are basis of $V^*$ of the form $(x^1,\dots,x^m,\xi^1,\dots,\xi^n)$, where
$(x^1,\dots,x^m)$ is a basis of $V_\0^*$, and $(\xi^1,\dots,\xi^n)$ a basis of
$V_\1^*$. Such a basis will be sometimes denoted by the symbol $(x,\xi)$. For
the corresponding predual basis $(e_1,\dots,e_m,f_1,\dots,f_n)$ of $V$, then
$(e_1,\dots,e_m)$ is a basis of $V_\0$ and $(f_1,\dots,f_n)$ is a basis of
$V_\1$. These notations will be used in particular for the canonical basis of
$\mathbb R^{(m,n)}$. Let $I= (i_1,
\dots,i_n)\in\{0,1\}^n$. Then we denote by $\xi^I$ the monomial
$(\xi^1)^{i_1}\dots (\xi^n)^{i_n}$ of $S(V^*)$. Let $\UC
\subset V_\0$ be an open set. Let $W$ be a Fr\'{e}chet supervector space. Then
any $ \phi\in \CC_V^\infty(\UC,W)$ is of the form
\begin{eqnarray}
\label{eq:funct}
\phi=\sum_I\xi^I \phi_I(x^1,\dots,x^m),
\end{eqnarray}
with $\phi_I$ is an ordinary $W$-valued smooth function defined in the
appropriate open subset of $\mathbb R^m$. Notice that $\phi_{\mathbb R}=
\phi_{(0,\dots,0)}(x^1,\dots,x^m)$   does not depend on the choice of the odd
coordinates $\xi^i$. We emphasize the fact that we write $\phi_I$ to the right
of  $\xi^I$ (recall that $\xi^I \phi_I= \pm \phi_I
\xi^I $, according to the sign rule).

We will denote by $\CC_{V,c}^\infty(\UC,W)$
  the subspace of $\CC_V^\infty(\UC,W)$ of function with compact support. Then
the  \emph{distributions on $V$ defined in $\UC$} are the elements of the
(Schwartz's) dual of $\CC_{V,c}^\infty(\UC )$. If $t$ is   a  distribution, we
will use the notation $$
t(\phi)=\int_{V } t(v) \phi(v) $$
for $\phi\in \CC_{V,c}^\infty(\UC)$. We will also use complex valued
distributions, defined in an obvious way.

A \emph{Berezin integral} (or  \emph{Haar}, or  \emph{Lebesgue}) is by
definition a distribution on $V$ which is  
invariant by translations (i.e.  which vanishes on functions of the form
$\de_X\phi$ where $\phi\in \CC^{\oo}_{V,c}(V)$ and $\de_X$ is the vector field
on $V$ with constant coefficients corresponding to $X\in V$: for $f\in V^*$,
$\de_{X}f=(-1)^{p(X)p(f)}f(X)$). 

Let $t$ be such a Berezin distribution. Let $\phi\in \CC^{\oo}_{V,c}(V)$. Let $\PC$ 
be a near superalgebra and  $a\in V_{\PC}$, then the function $\phi_{a}(v)=\phi(v+a)$ 
is a well defined function $\phi_{a}\in\CC^{\oo}_{V,c}(V,\PC)$. Then,
Taylor formula (\ref{FonctLisse}) and invariance by translation 
of $t$ implies that:
\begin{equation}
\int_{V}t(v)\phi_{a}(v)=\int_{V}t(v)\phi(v+a)=\int_{V}t(v)\phi(v).
\end{equation}

Up to a multiplicative constant, there is
exactly one    Berezin integral, and it is an important matter in this article
to choose a particular one for the symplectic oriented supervector spaces (see
below).

\medskip A choice of a standard system of coordinates determines a specific
choice $d_{(x,\xi)}$ of a Berezin integral by the formula
\begin{equation}
\label{eq:haar}
\begin{split}
 \int_V d_{(x,\xi)}(v)\phi(v)&=(-1)^{\frac{n(n-1)}2}\int_{\mathbb R^m}
   \big|\,dx^1\ldots dx^m\big| \phi_{(1,\dots,1)}(x^1,\dots,x^n)\\
&=\int_{\mathbb{R}^m}  \big|\,dx^1\ldots dx^m\big|
\Big(\frac{\de}{\de \xi^1}\dots\frac{\de}{\de \xi^n}
\phi\Big)_\mathbb R (x^1,\dots,x^m),
\end{split}
\end{equation}
for $\phi\in \CC_{V,c}^\infty(V_\0)$, where $\big|\,dx^1\ldots dx^m\big|$ is the
Lebesgue measure on $\mathbb{R}^m$. Note that this formula can also be applied
to any $\phi\in
\CC_{V,c}^\infty(\UC,W)$ with a result in $W$.

The choice of sign is such that Fubini's formula holds. More precisely, let
$V,W$ be two supervector spaces of dimensions $(m,n)$ and $(p,q)$. Let
$(x^1,\dots,x^m,\xi^1,\dots,\xi^n)$ be standard coordinates on $V$ and
$(y^1,\dots,y^p,\eta^1,\dots,\eta^q)$ be standard coordinates on $W$. Then
$(x,y,\xi,\eta)=(x^1,\dots,x^m,y^1,\dots,y^p,\xi^1,\dots,\xi^n,\eta^1,\dots%
,\eta^q)$ defines standard coordinates on $V\times W$. Let $\phi(v,w)$ is a
smooth compactly supported function on $V\times W$. Then:
\begin{equation}
\label{eq:Fubini}
\int_{V\times W}d_{(x,y,\xi,\eta)}(v,w)\phi(v,w)=\int_Vd_{(x,\xi)}(v)\Big(\int_W
d_{(y,\eta)}(w)\phi(v,w)\Big).
\end{equation}
We write:
\begin{equation}
\label{eq:Fubini-int} d_{(x,y,\xi,\eta)}(v,w)=d_{(x,\xi)}(v)d_{(y,\eta)}(w).
\end{equation}

In particular, since $V=V_\0\oplus V_\1$, $d_{x,\xi}=d_x d_\xi$ and formula
(\ref{eq:haar}) is a particular case of formula (\ref{eq:Fubini}).

\medskip Let us stress that if $V_\1$ is not $\{0\}$, in the setting of
supermanifolds  there is no natural notion of measure  on $V$ and no natural
notion of positive distribution on $V$. Thus we use these notions  only for
(ungraded) vector spaces, or for the even part $V_\0$ of a supervector space,
which is then regarded as an  ungraded  vector space. Otherwise, we use the
terms \emph{distribution} or \emph{integral}.

\medskip In this article, we will be in fact interested by complex valued
distributions. Then we allow   standard basis $(e_1,\dots,e_m,f_1,\dots,f_n)$ of
$V\otimes \mathbb C$, where $(e_1,\dots,e_m)$ is a basis of $V_\0$ and
$(f_1,\dots,f_n)$ is a basis of $V_\1\otimes \mathbb C$. Then the  dual basis
$(x,\xi)$ provides a coordinate system $(x)$ on $V_\0$ and a dual basis $(\xi)$
of $V_\1^*\otimes \mathbb C$. Any $f \in
\CC_{V,c}^\infty(V_\0,\mathbb C)$ can be written in the form (\ref{eq:funct}),
and the (complex) Berezin integral $d_{(x,\xi)}$ is again well defined by  
formula (\ref{eq:haar}).

\subsection{Generalized functions}
\subsubsection{Definition} Let $V$ be a finite dimensional supervector space and
$\UC\subset V_0$ be an open set. Let $(x,\xi)$ be a standard coordinates system
on $V$.  As usual, we will say that a distribution $t$ on $V$ defined in $\UC$
is \emph{smooth} (resp. \emph{smooth compactly supported}) if there is a
function $\psi\in\CC^{\oo}_V(\UC)$ (resp. $\psi\in\CC^{\oo}_{V,c}(\UC)$) such
that $t(v)=d_{x,\xi}(v)\psi(v)$. It means that for any
$\phi\in\CC^{\oo}_V(\UC)$:
\begin{equation}
t(\phi)=\int_Vd_{x,\xi}(v)\psi(v)\phi(v).
\end{equation}
This definition does not depend on the standard coordinates system $(x,\xi)$.

\medskip

By definition, the \emph{generalized functions on $V$ defined on $\UC$} are the
elements of the (Schwartz's) dual of the space of smooth compactly supported
distributions. For a generalized function $\phi$ and a smooth compactly
supported distribution $t$, we write:
\begin{equation}
\phi(t)=(-1)^{p(t)p(\phi)}\int_V t(v)\phi(v).
\end{equation}
(the spaces of distributions and thus of generalized functions are naturally
$\mathbb{Z}/2\mathbb{Z}$-graded.)

We denote by $\CC^{-\oo}(\UC)$ the space of generalized functions on  $\UC$ and
by $\CC^{-\oo}_V(\UC)$ the space of generalized functions on $V$ defined on
$\UC$.

Let us remark that, as $\CC^\oo_V(\UC)=\CC^{\oo}(\UC)\tens\Lambda( V_\1^*)$, we
have:
\begin{equation}
\CC^{-\oo}_V(\UC)=\CC^{-\oo}(\UC)\tens\Lambda (V_\1^*).
\end{equation}

\medskip

Let $W$ be a Fr\'echet supervector space. A $W$-valued generalized function  is
a continuous homomorphism (in sense of Schwartz) from the space of smooth
compactly supported distributions to $W.$ We denote by $\CC^{-\oo}_V(\UC,W)$ the
set of $W$-valued generalized functions. If $W$ is finite dimensional, we have
$\CC^{-\oo}_V(\UC,W)=\CC^{-\oo}_V(\UC)\tens W$. We will be  in particular
concerned with the cases $W=\mathbb{C}$ and $W=\L(E^*)$ for some finite
dimensional vector space $E$.

\medskip

\subsubsection{Wave front set}

\begin{defi} Let $V$ be a supervector space. Let $\UC\subset V_\0$ be an open
subset of $V_\0$.  Let $\psi\in\CC_V^{-\oo}(\UC)$ be a generalized function on
$V$ defined on $\UC$. Let $(\xi^1,\dots,\xi^n)$ be a basis of $V_\1^*$. We put
$\psi=\som_{I}\xi^I\psi_I$ where $\psi_I$ is a generalized function on $\UC.$ 

The wave front set  $WF(\psi)$ of $\psi$, is by definition the union of the wave
front sets of the $\psi_I$:
\begin{equation}
WF(\psi)=\union_{I}WF(\psi_I)\subset T^*\UC,
\end{equation}
where $WF(\psi_I)$ is the wave front set of $\psi_I$ and $ T^*\UC=\UC\times
V_\0^*$ is the cotangent bundle of $\UC$.
\end{defi}

The definition of $WF(\psi)$ does not depends on the choice the basis
$(\xi^{i})$ of  $V_{\1}^*$.

\subsection{Rapidly decreasing functions}

 We say (cf. for example \cite[Chapter 7]{Hor83}) that
$\phi\in\CC^{\oo}(\mathbb{R}^m)$ is rapidly decreasing if for any
$(\a_1,\dots,\a_m)\in\mathbb{N}^m$ and any $(\be_1,\dots,\be_m)\in\mathbb{N}^m$,
\begin{equation}
\sup\Big|(x^1)^{\be_1}\dots(x^m)^{\be_m}
\frac{\de^{\a_1}}{\de (x^1)^{\a_1}}\dots \frac{\de^{\a_m}}{\de
(x^m)^{\a_m}}\phi\Big|<+\oo.
\end{equation}
where $(x^1,\dots,x^m)$ are the canonical coordinates on $\mathbb{R}^m$.

\begin{defi} Let $V$ be a supervector space. Let $\phi\in\CC_V^{\oo}(V_{\0})$ be
a smooth function on $V$. Let $(x^1,\dots,x^m,\xi^1,\dots,\xi^n)$ be a basis of
$V^*$. We put $\phi=\som_{I}\xi^I\phi_I(x^1,\dots,x^m)$ where
$\phi_I\in\CC^\oo(\mathbb{R}^m)$.

We say that $\phi$ is rapidly decreasing if for any $I$, $\phi_I$
 is a rapidly decreasing function on $\mathbb{R}^m$.
\end{defi} 

This definition does not depend on the choice of the basis $(x^{i},\xi^j)$ of
$V^*$. 

\subsection{Holomorphic functions} Let $V$ be a complex finite dimensional
supervector spa\-ce.  Let $\UC$ be an open subset of $V_{\0}$. Let $W$ be a
complex Fr\'echet supervector space. By definition, an holomorphic function on
$V$ with values in $W$ defined on $\UC$ is an holomorphic function on $\UC$
with values in $\L(V_{\1}^*)\tens W$. We denote by $\HC_{V}(\UC,W)$ (resp.
$\HC(V,W)$) the algebra of holomorphic functions on $V$ (resp. on $V_{\0}$)
with values in $W$ defined on $\UC$:

\begin{equation}
\HC_{V}(\UC,W)=\HC(\UC,W)\tens \L(V_{\1}^*).
\end{equation}

Let $\PC$ be a (real) near superalgebra. We denote by $\PC_{\mathbb
C}=\PC\tens\mathbb C$ its complexification. We put $V_{\PC}=(V\tens \PC_{\mathbb
C})_{\0}$ and $W_{\PC}=(W\tens \PC_{\mathbb C})_{\0}$. As in the real case we have a 
body map $\bb:V_{\PC}\to V_{\0}$ which extends the canonical projection 
$\bb\tens 1:\PC\tens \mathbb C\to\mathbb R\tens \mathbb C\simeq\mathbb C$. 
We put for $\UC$ open in $V_{\0}$:
 $V_{\PC}(\UC)=\{v\in V_{\PC}\,/\, \bb(v)\in \UC\}$.

Let $\phi\in\HC_{V}(\UC,W)$ and $\a\in V_{\PC}(\UC)$ as in the real
case, we denote by $\phi(\a)\in W_{\PC}$ the image of $\a$ by $\phi$ defined by
a formula analogous to formula (\ref{FonctLisse}). The map $\a\mapsto\phi(\a)$
defines an holomorphic function $\phi_{\PC}\in\HC(V_{\PC},W_{\PC})$ on $V_{\PC}
$ with values in $W\tens{\PC_{\mathbb C}}$.

As for smooth functions on a real supervector space, if $\PC$ is  large enough,
the map $\phi\mapsto\phi_{\PC}$ is injective.

\subsection{Analytic functions}

Let $V$ be a real finite dimensional supervector space. Let $\UC$ be an open subset of $V_{\0}$.
Let $W$ be a real  supervector space.
We denote by $\CC^{\w}(\UC,W)$ the set of analytic functions on $\UC$ with values in $W$.
 We put:
\begin{displaymath}
\CC^{\w}_{V}(\UC,W)=\CC^{\w}(\UC,W)\tens\L(V_{\1}^*).
\end{displaymath}
We call the elements of $\CC^{\w}_{V}(\UC,W)\subset\CC^{\oo}_{V}(\UC,W)$ 
the analytic functions on $V$
 with values in $W$ defined on $\UC$.

As usual we put $\CC^{\w}(V,W)=\CC^{\w}_{V}(V_{\0},W)$ and 
$\CC^{\w}_{V}(\UC)=\CC^{\w}_{V}(\UC,\mathbb R)$ (resp. 
$\CC^{\w}(V)=\CC^{\w}(V,\mathbb R)$).

\subsection{Supertrace and Berezinians}\label{sec:Ber} Let $V$ be a supervector
space. We denote by $\g{gl}(V)$ the Lie superalgebra of endomorphisms of $V$.

Let $\AC$ be a commutative superalgebra. We write an element of
$\g{gl}(V)_{\AC}$ in the form
\begin{equation}
M=
\begin{pmatrix} A&B\\
C&D\\
\end{pmatrix}\in\g{gl}(V)_{\AC}.
\end{equation}
where $A\in \g{gl}(V_0)\otimes \AC_0$, $D\in \g{gl}(V_1)\otimes
\AC_0$, $B\in Hom(V_\1,V_\0)\otimes \AC_\1$, and  $C\in Hom(V_\0,V_\1)\otimes
\AC_\1$

We recall the definition of the supertrace:
\begin{defi} The supertrace of $M\in\g{gl}(V)_{\AC}$ is defined by
\begin{equation}
\str(M)=\tr(A)-\tr(D),
\end{equation}
where $\tr$ is the ordinary trace.
\end{defi}

\medskip

Berezin introduced the following generalizations of the determinant (cf.
\cite{Ber87,BL75,Man88}), called the
\emph{Berezinian} and \emph{inverse Berezinian}.

\begin{defi} If $D$ is invertible, we define:
\begin{equation}
\label{eq:bera}
\Ber(M)=\det(A-BD^{-1}C)\det(D)^{-1},
\end{equation}
and if $A$ is invertible:
\begin{equation}
\label{eq:beram}
\Ber^-(M)=\det(A)^{-1}\det(D-CA^{-1}B).
\end{equation}
\end{defi}

\begin{defi} Assume moreover that  $\AC$ is a near superalgebra. If $D$ is
invertible, we define (cf. \cite{Vor91}):
\begin{align}\label{eq:berd}
\Ber_{(1,0)}(M) &=\Big|\det(A-BD^{-1}C)\Big|\det(D)^{-1},
\end{align} and if $A$ is invertible:
\begin{align}\label{eq:berdm}
\Ber_{(1,0)}^-(M) &=\Big|\det(A)^{-1}\Big|\det(D-CA^{-1}B),
\end{align}
\end{defi}

All these functions are multiplicative, and when  both $A$ and $D$ are
invertible, it is known that $\Ber^-(M)=\Ber^{-1}(M)$ and
$\Ber_{(1,0)}^-(M)=\Ber_{(1,0)}^{-1}(M)$.

Recall that $\g{gl}(V)_\0$ consists of the matrices $
\begin{pmatrix}A&0\\
0&D
\end{pmatrix}$ with $A\in \g{gl}(V_\0)$ and $D\in \g{gl}(V_\1)$. We consider the
two open sets $\UC'=GL(V_\0)\times
\g{gl}(V_\1)$, and $\UC"= \g{gl}(V_\0)\times GL(V_\1)$. Formula (\ref{eq:bera})
defines a rational function on the open set $\UC"$  of the supermanifold
$\g{gl}(V)$. Formula (\ref{eq:berd}) defines a smooth function on the open set
$\UC"$  
of the supermanifold $\g{gl}(V)$. We still denote by $\Ber$ and $\Ber_{(1,0)}$  
the elements of $\CC^\infty_{\g{gl}(V)}(\UC")$ whose evaluation in
$\g{gl}(V)_\AC$ is given as above. We similarly define the elements 
$\Ber^-$ and $\Ber_{(1,0)}^-$     of $\CC^\infty_{\g{gl}(V)}(\UC')$.

\subsection{Symplectic oriented supervector spaces} \label{Symp} Let
$V=V_\0\oplus V_\1$ be a supervector space. A
\emph{symplectic form} $B$ on $V$ is a non degenerate even skew symmetric
bilinear form on $V$. It means that $V_\0$ and $V_\1$ are orthogonal, that the
restriction of $B$ to $V_\0$ is a non degenerate   skew symmetric bilinear form,
and that the restriction of $B$ to $V_\1$ is a non degenerate     symmetric
bilinear form. We call $V$, provided with $B$, a
\emph{symplectic supervector space}.

\medskip Such a space is a direct sum of $(2,0)$-dimensional symplectic
supervector spaces (i.e. $2$-dimensional symplectic  vector spaces), and of
$(0,1)$-dimensional   symplectic supervector spaces  (i.e. $1$-dimensional
quadratic 
vector spaces). We first review these building blocks.

\subsubsection{Symplectic $2$-dimensional vector spaces}\label{2symp} Let
$V=V_\0$ a purely even $2$-dimensional  
symplectic space. A
\emph{symplectic basis} of $V$ is a basis $(e_1,e_2)$ such that $B(e_1,e_2)=1$,
$B(e_1,e_1)=0$,  $B(e_2,e_2)=0$.
 It defines a dual \emph{symplectic coordinate system} $(x^1,x^2)$, an
orientation of $V^*$, and a
\emph{Liouville integral} (a particular normalization of the Berezin integral):
$$
\phi\in \CC^\infty_c(V)=\CC^\infty_c(\mathbb R^2)\mapsto
 \int_V d_V(v) \phi(v)
 =\frac{1}{2\pi }\int |dx^1  dx^2|\phi(x^1,x^2).
 $$

\subsubsection{Symplectic $1$-dimensional odd vector spaces} Let $V=V_\1$ a
purely odd $1$-dimensio\-nal   symplectic supervector space (i.e. a
$1$-dimensional quadratic space). A symplectic basis of $V$ is a basis $(f  )$
such that $B(f ,f )= 1$. However, such a basis does not always exists, and we 
allow $f$ to be in $V_\1
\cup \ii V_\1 \subset V_\1\otimes \mathbb C$. Let $(\xi)\in V_\1^*
\cup \ii V_\1^* $ be the dual basis. It defines a \emph{Liouville integral}
(which is complex valued if $B$ is negative definite) $d_{V }$:
 $$
\phi=a+\xi b\in \Lambda(V^*\otimes \mathbb C)  \mapsto \int_V d_{V}(v)
\phi(v)=b.$$

We will call the choice of $(\xi)$ (the other possible choice is $(-\xi)$) an
\emph{orientation} of $V_\1$. If $B$ is positive definite, then $(\xi)$ is a
basis of $V_\1^*$, and so defines an orientation in the usual sense.  If $B$ is
negative definite, then $(-\ii\xi)$ is a basis of $V_\1$, and so defines an
orientation in the usual sense.

\subsubsection{General case} Let us go back to the general case. Since $V_\0$ is
a classical symplectic space, there is a canonical normalization of Lebesgue
integral on $V_\0$, the Liouville integral,   which we recall now.

The dimension $m$ of $V_\0$ is even. We choose a  symplectic basis
$(e_1,\dots,e_m)$ of $V_\0$, that is  $V_\0$ is the direct sum of $m/2$
symplectic vector spaces generated by the pairs $(e_1,e_2) ,
 (e_3,e_4) , \dots$, and  $B(e_1,e_2)=1, B(e_3,e_4)=1, \dots$. The dual basis
$(x^i)$ of $V_\0^*$ is called a symplectic coordinate system. The Liouville
integral on $V_\0$ is $$
\frac{1}{(2
\pi)^{m/2}}|dx^1\dots dx^m|.$$
  The Liouville integral does not depend on the choice of the symplectic basis
of $V_0$.

\medskip On $V_\1$, we define a symplectic basis $(f_1,\dots f_n)$ as an
orthonormal basis of $V_\1\otimes \mathbb C$ such that $f_i\in V_\1$ or $f_i\in
\ii V_\1$ for all $i$. Let $(\xi^1,\dots\xi^n)$ be the dual basis. The pair of
functions $\pm \xi^1\dots\xi^n $ does not depend on the symplectic basis
$(f_1,\dots,f_n)$. A choice of one of the two elements of $\pm \xi^1\dots\xi^n $
is called an \emph{orientation} of $V_\1$. If $V_\1$ is oriented, an oriented
symplectic coordinate system on $V_\1$  is a basis for which the orientation is
$\xi^1\dots\xi^n $.

We call the corresponding  Berezin integral $d_{\xi}$ the Liouville integral of
the oriented symplectic space $V_\1$.

\medskip Let us   remark that in the specially interesting case where $B$ is
positive definite,   a symplectic basis is a basis of $V_\1$, and not only of
$V_\1\otimes \mathbb C$.

\medskip We define an \emph{oriented symplectic supervector space} as a
symplectic supervector space $(V,B)$ provided with an orientation of $V_\1$.   A
symplectic oriented basis $(e_1,\dots,e_m,f_1,\dots,$ $f_n)$ is a basis of
$V\otimes
\mathbb C$ such that $(e_1,\dots,e_m)$ is a symplectic basis of $V_\0$, and
$(f_1,\dots,f_n)$ an oriented symplectic basis of $V_\1\otimes
\mathbb C$. We use the corresponding dual system of coordinates $(x,\xi)$, and
the associated Berezin integral $d_{(x,\xi)}$ will be denoted by
$d_V$. Then $V$ becomes a symbol   bearing a supermanifold  structure, a
symplectic structure, an orientation...

\subsection{Symplectic Lie superalgebras} Let $V=(V,B)$  be a symplectic
supervector space. We denote by $\g{spo}(V,B)$ (or $\g{spo}(V)$) the Lie
subsuperalgebra of $\g{gl}(V)$ consisting of  endomorphisms of $V$ which leave
$B$ invariant. We call it the \emph{symplectic Lie superalgebra} ---it is
usually called \emph{orthosymplectic}---.

We have $\g{spo}(V,B)_\0=\g{sp}(V_\0)\oplus\g{so}(V_\1)$.

A particular  linear even isomorphism $\mu$ of $\g{spo}(V)$ to $ S^2(V^*)$,
called the moment mapping,  will play an important role. Thus $\mu$ is an
element of $S(\g{spo}(V)^*\tens V^*)$, and we consider it as a function   on the
supermanifold $\g{spo}(V)\times V$, linear in the first variable, polynomial of
degree $2$ in the second variable, and globally homogeneous of degree $3$. It is
defined by the formula $$
\mu(X,v)=-\frac{1}{2} B(v,Xv), $$
where, for any commutative superalgebra $\AC$,  $X$ and $v$ are $\AC$-valued
points of $\g{spo}(V)$ and $V$, and $B(v,Xv)\in \AC$ is defined by the natural
extension of scalars. Considering a basis $G_k$ of $\g{spo}(V)$, the dual basis
$Z^k$, the generic point $X=G_k Z^k$, a basis $g_i$ of $ V $, the dual basis
$z^i$, and the generic point $v=g_i z^i$, we obtain: $$
\mu =-\frac{1}{2} B(g_i,G_k g_j)z^j  Z^k z^i. $$

We will also consider $e^\mu$, which is a smooth (and even analytic) function on
the supermanifold $\g{spo}(V)\times V$.

\medskip

Let us explain the choice of the constant $-\frac12$ in definition of $\mu$ and
why we call $\mu$ the moment mapping. 

The symplectic form on $V$ gives to the associated supermanifold a structure of
symplectic supermanifold. 
We define a Poisson bracket on $S(V^*)$ by the following. Let $f\in V^*$, we
denote by $v_{f}$ the element of $V$ such that for any $w\in V$:
\begin{equation}
\label{eq:IsoSymp} 
B(v_{f},w)=f(w).
\end{equation}
This gives an isomorphism from $V^*$ onto $V$. For $f,g\in V^*$, we put:
\begin{equation}
\{f,g\}=B(v_{f},v_{g})
\end{equation}
and we extend it to a Poisson bracket on $S(V^*)$.

Let $\check\mu$ be the linear form on $\g g$ with values in $S^2(V^*)$ such that
$\check\mu(X)(v)=\mu(X,v)$.

With the above definitions we have:
\begin{equation}
\label{eq:Poisson}
\{\check\mu(X),\check\mu(Y)\}=\check\mu([X,Y]).
\end{equation}
Thus $\check\mu$ is a morphism of Lie algebras.

\bigskip

We extend $\check \mu$ to a morphism of superalgebras from $S(\g{spo}(V))$ to
$S(V^*)$. 
 More precisely, we put for $X_1\dots X_k\in S^k(\g{spo}(V))$:
\begin{equation}
\label{eq:checkmu}
\check\mu(X_1\dots X_k)=\check\mu(X_1)\dots\check\mu(X_k).
\end{equation}

Its image is $\sdir_{k\in\mathbb N}S^{2k}(V_{\0}^*)$.
 Indeed, the natural morphism $S^k(S^{2}(V^*))\to S^{2k}(V^*)$ is surjective.
Moreover, $\check\mu:\g{spo}(V)\to S^2(V^*)$ is bijective and by definition
$\check\mu$ factorises:
\begin{equation}
S^k(\g{spo}(V))\to S^k (S^2(V^*))\to S^{2k}(V^*).
\end{equation}

\medskip We put $\IC_{k}=\ker(\check\mu)\cap S^k(\g{spo}(V))$. We have:
$\ker(\check\mu)=\sdir_{k\in\mathbb N}\IC_{k}$. We choose a supplementary
$\SC_{k}\subset S^k(\g{spo}(V))$ of $\IC_{k}$:
\begin{equation}
\label{eq:SCk}
S^k(\g{spo}(V))=\IC_{k}\oplus\SC_{k}.
\end{equation}
 Let 
 \begin{equation}
 \label{eq:SC}
\SC=\sdir_{k\in\mathbb N}\SC_{k},
\end{equation}
we have $S(\g{spo}(V))=\SC\oplus\ker(\check\mu)$. 

Then the restriction $\check\mu:\SC\to \sdir_{k\in\mathbb N}S^{2k}(V^*)$ is
bijective. We denote by 
\begin{displaymath}\label{Xi}
\Xi: \sdir_{k\in\mathbb N}S^{2k}(V^*)\to\SC\subset S(\g{spo}(V))
\end{displaymath}
 its inverse. For $P\in S^{2k}(V^*)$, $\Xi(P)\in\SC_{k}\subset S^k(\g{spo}(V))$
and
\begin{equation}
 \check\mu(\Xi(P))=P.
\end{equation}

\medskip

\section{SuperPfaffian I : an analytic function}\label{analytic}

Let $V=V_\0\oplus V_\1$ be an oriented  supervector symplectic space of
dimension $(m,n)$. We already defined the symplectic integral $d_V$ and the
moment map $\mu$. We now define some open subsets of $\g{spo}(V)_\0 $ in which
we want to define functions on $\g{spo}(V)$.

For $X\in \g{sp}(V_\0 )$, $v \mapsto B(v,Xv)$ is a quadratic form on $V_\0$. We
denote by $\UC\subset \g{sp}(V_\0 )$ the open set of $X
\in \g{sp}(V_\0 )$ for which this form is non degenerate (or equivalently,
$X\in\g{sp}(V_\0 )\subset\g{gl}(V_\0)$ is invertible). It is the disjoint union
of the subsets $\UC_{p,q}$ (with $p+q=m$) where $(p,q)$ is the signature of the
quadratic form. Let us also use the notation $\UC^+= \UC_{m,0}$ (resp. $\UC^-=
\UC_{0,m}$) for the open set of $X$ for which it is positive definite (resp.
negative definite). It is an open convex cone in $\g{sp}(V_\0 )$. We denote by
$\VC=\UC\times \g{so}(V_\1 )$, $\VC_{p,q}=\UC_{p,q}\times \g{so}(V_\1 )$,
$\VC^+=\UC^+\times
\g{so}(V_\1 )$,  $\VC^-=\UC^-\times
\g{so}(V_\1 )$ the corresponding open subsets of $\g{spo}(V)_\0 $. Then $\VC^+$
is an open convex cone in $\g{spo}(V)_\0 $.

In this section, we prove:

\begin{theo}\label{th:Spf} There exists a unique function 
$$
\Spf  \in \CC^\w_{\g{spo}(V)}(\VC^+,\mathbb C) $$
such that for any near superalgebra $\PC$ and any  element $X$  of
$\g{spo}(V)_\PC$ such that the body $\bb(X)$ is in $\VC^+$, (that is, with the
notations of section \ref{superman}, $X\in\g{spo}(V)_\PC(\VC^+)$):

\begin{equation}
\label{eq:Spf}
\Spf (X)=\int_{V}d_V(v)\exp(\mu(X,v))=\int_{V}d_V(v)\exp(-\frac{1}{2} B(v,Xv)).
\end{equation}
\end{theo}

We will prove it theorem in section \ref{correct}. We call the function $\Spf$ 
defined in the preceding theorem the
\emph{superPfaffian}. We begin by giving typical examples.

\subsubsection{Symplectic $2$-dimensional vector spaces}\label{2symp2} Let
$V=V_\0$ a purely even $2$-dimensional  
symplectic space (cf. subsection \ref{2symp}). Choose  a  symplectic basis 
$(e_1,e_2)$ of $V$.  The  dual coordinate system  will be denoted by $(x,y)$
(instead of $(x^1,x^2)$). An element $X\in\g{spo}(V)$ is represented by a matrix
$X=
\begin{pmatrix}
      a &b    \\
       c&  -a
\end{pmatrix}\in\g{sl}(2,\mathbb R)$. The corresponding quadratic form $B(v,Xv)$
is equal to $B(e_1x +e_2y, X(e_1x +e_2y))= B(e_1x +e_2y,( e_1 a + e_2 c)x +(e_1
b - e_2 a)y))= x(cx -a y)-y(a x + b y)= c x^2 -2a x  y -by^2$. Thus, if
$v=e_{1}x+e_{2}y$ is the generic point of $V$:
\begin{equation}
\label{eq:mu2symp}
\mu(X,v)=\frac{-1}2(c x^2 -2a x  y -by^2).
\end{equation}

 The set $\UC^+$ is defined by the equations $\det(X)=-a^2-bc>0$, and $c>0$. It
is one of the two connected components of the set of invertible elliptic
matrices. For example, the matrix $\begin{pmatrix}
      0 &-c    \\
       c&  0
\end{pmatrix}$ is in $\UC^+$ if $c>0$, and these matrices are representatives
for the conjugacy classes of $SL(2,\mathbb R)$ in $\UC^+$.

For $X\in \UC^+$ the integral
\begin{equation}
\Spf (X)=\int_V d_V(v) \exp(-\frac{1}{2} B(v,Xv))
\end{equation}
is convergent and defines an analytic function on $\UC^+$. By a suitable change
of  variable  we see that the function $\Spf $ is invariant by $SL(2,\mathbb
R)$. For $X=
\begin{pmatrix}
      0 &-c    \\
       c&  0
\end{pmatrix}$ with $c>0$, we get
\begin{equation*}
\Spf (X)=\frac{1}{2 \pi}\int | dx\, dy| \exp(-\frac{1}{2} c(x^2+y^2)),
\end{equation*} and so
\begin{equation}
\Spf(
\begin{pmatrix}
      0 &-c    \\
       c&  0
\end{pmatrix})= 1/c \quad \text{if  $c>0$}.
\end{equation}
Notice also that $\Ber^-(X)=1/\det(X)=1/c^2$.

\medskip

Our conclusion is that  \emph{on $\UC^+$ the inverse Berezinian is positive, and
that $\Spf $ is the positive square root of the inverse Berezinian}.

\medskip Consider now a near superalgebra $\PC$. Let  $X\in
\g{sl}(2,\mathbb R)\otimes \PC_\0$   such that $\bb(X)\in \UC^+$. We write 
$X=\begin{pmatrix}
      a+\alpha & b+\beta    \\
       c+\gamma&  -a-\alpha
\end{pmatrix}$, with $a,b,c$ real such that $\det(X)=-a^2-bc>0$ and $c>0$, and
with $\alpha,\beta,\gamma $ nilpotent elements of $\PC_\0$. Then, for
$v=e_1x+e_2y$, we obtain 
\begin{equation}
\begin{split}
\exp(-\frac{1}{2} B(v,Xv))&= \exp(-\frac12
((c+\gamma)x^2-2(a+\a)xy-(b+\be)y^2)\\
&=\exp(-\frac12 (cx^2-2axy-by^2))\exp(-\frac12 (\gamma x^2-2\a xy-\be y^2)).
\end{split}
\end{equation}
Since $\a,\be$ and $\gamma$ are nilpotent the last exponential is a polynomial
on $V$ with values in $\PC_0$. It follows
 that $\int_Vd_V(v)\exp(-\frac{1}{2} B(v,Xv))$ converges because
$cx^2-2axy-by^2$ is a positive definite quadratic form on $V$.

Now, on any compact set included in $\UC^+$,  $\exp(-\frac12 (cx^2-2axy-by^2))$ 
and its derivatives can be uniformly bounded by a rapidly decreasing  function
on $V$. Thus,  $\Spf$ is a smooth function  on $\UC^+$ and its derivatives are
obtained by derivation under the summation symbol. It follows that the value of
$\Spf\begin{pmatrix}
      a+\alpha & b+\beta    \\
       c+\gamma&  -a-\alpha
\end{pmatrix}$ obtained by applying formula (\ref{FonctLisse}) coincides we the
value obtained by integration of formula (\ref{eq:Spf}). In particular it
implies  that equality $\Spf(X)^2=\Ber(X)$ is still valid on
$\g{sp}(V)_\PC(\UC^+)$ (that is for $X=\begin{pmatrix}
      a+\alpha & b+\beta    \\
       c+\gamma&  -a-\alpha
\end{pmatrix}$ with $a,b,c$ real such that $\det(X)=-a^2-bc>0$ and $c>0$, and
with $\alpha,\beta,\gamma $ nilpotent elements of $\PC_\0$).

\subsubsection{Symplectic oriented $2$-dimensional odd vector
spaces}\label{sec:Ex2Odd} If $V$ is of dimension $(0,1)$, then  $\g{spo}(V)$ is
$\{0\}$, so the next interesting example is when $V$ is of dimension $(0,2)$.

Then $V=V_\1$, and we choose an oriented symplectic basis $(f_1,f_2)$ of
$V_\1\otimes \mathbb C$. Recall that it means in particular that
$B(f_1,f_2)=0,B(f_1,f_1)=B(f_2,f_2)=1$. We denote the dual basis by
$(\xi,\eta)$. Then the orientation is the choice of $\xi\eta\in
S^2(V^*)=\Lambda^2 (V_\1^*)$ (versus $-\xi\eta$). In this case the superPfaffian
is the ordinary Pfaffian, a polynomial (and in fact, a linear) function on
$\g{so}(V)$ and so it is also defined on $\g{so}(V\otimes \mathbb C)$. An
element $X$ of $\g{so}(V\otimes \mathbb C)$ is represented in the given basis by
a matrix $X=\begin{pmatrix}
      0 &-c    \\
       c&  0
\end{pmatrix}$ with $c \in \mathbb C$. Then, for the generic point
$v=f_1\xi+f_2\eta$, we have $B(v,Xv)= B(f_1\xi+f_2\eta,f_2 c\,\xi-f_1 c \,\eta)= 2 c
\xi \eta$ and $\exp(-\frac{1}{2}B(v,Xv))=1-c\xi \eta$.

We obtain that  the integral $\int_V d_{(\xi,\eta)}
(v)\exp(-\frac{1}{2}B(v,Xv))$, which by definition is the constant term of the
function $\frac{
\partial}{
\partial\xi}\frac{
\partial}{
\partial\eta} (1-c\xi \eta)$, is equal to $c$, and so
\begin{equation} 
\Spf (
\begin{pmatrix}
      0 &-c    \\
       c&  0
\end{pmatrix})= c \quad \text{if  $c\in \mathbb C$}.
\end{equation}
Notice also that $\Ber^-(X)= \det(X)= c^2$. So again $\Spf (X)$ is a square root
of the inverse Berezinian.

\medskip

If $B$ is positive (resp. negative) definite, $f_1,\,f_2\in V_\1$ (resp. $\in\ii
V_\1$) and thus if $X\in\g{so}(V_1)$, its matrix in the basis $(f_1,f_2)$ is
real and $\Spf(X)\in\mathbb{R}$.

If $B$ is hyperbolic, we can take $f_1\in V_\1$ and $f_2\in \ii V_\1$ and thus
if $X\in\g{so}(V_\1)$, its matrix in the basis $(f_1,f_2)$ is purely imaginary
and $\Spf(X)\in\ii\mathbb{R}$.

\subsection{Proof of theorem \ref{th:Spf}}\label{correct} 

Let $\PC$ be any near superalgebra. Then $X\in\g{spo}(V)_\PC(\VC^+)$  means
$X\in \g{spo}(V)_\PC$ with $\bb(X)\in\VC^+$.

We denote by $\Spf_{\PC}$ the function on $\g{spo}(V)_\PC(\VC^+)$ such that for
any $X\in \g{spo}(V)_\PC(\VC^+)$, 
\begin{displaymath}
\Spf_{\PC}(X)=\int_{V}d_V(v)\exp(\mu(X,v))=\int_{V}d_V(v)\exp(-\frac{1}{2}
B(v,Xv)).
\end{displaymath}

\subsubsection{$\Spf_\PC$ is a well defined $\PC$-valued analytic function on
$\g{spo}(V)_\PC(\VC^+)$}
\label{Spf_P}
 Let $X\in\g{spo}(V)_\PC(\VC^+)$. Let  $X_0\in\UC^+$ and
$X_1\in\g{sp}(V_\1)$ such that $\bb(X)=X_0+X_1$. Thus $X=X_{0}+X_{1}+N$ with $N$
nilpotent.

Let $(e_1,\dots,e_m,f_1,\dots,f_n)$ be a standard basis of $V$. Let $(x,\xi)$ be
the dual basis. Let $v=e_ix^i+f_j\xi^j$ be the generic point of $V$, 
$v_\0=e_ix^i$ be the generic point of $V_\0$ and $v_\1=f_j\xi^j$ the generic
point of $V_\1$. 
We have $v=v_\0+v_\1$, $B(v,X_0v)=B(v_\0,X_0v_\0)$ and
$B(v,X_1v)=B(v_\1,X_1v_\1)$. 

In particular $B(v,X_1v)$ is nilpotent. Thus $B(v,(X-X_0)v)\in
(S^2(V^*)\tens\PC)_\0$ is nilpotent. It follows that $\exp(-\frac12
B(v,(X-X_0)v))\in (S(V^*)\tens\PC)_\0$ is a polynomial on $V$ with values in
$\PC$ and that $X\mapsto\exp(-\frac12 B(v,(X-X_0)v))\in (S(V^*)\tens\PC)_\0$ is
polynomial on $\g{spo}(V)_\PC$.

Let $Z\in \g{spo}(V)_{\PC}$ such that $\bb(Z)\in\g{so}(V_{\1})$. Then $B(v,Zv)$
is a nilpotent element of $S^2(V^{*})_{\PC}$. We put:
\begin{equation}
\label{P} P(Z,v_\0)=
\int_{V_\1}d_{V_\1}(v_\1)\exp(-\frac12 B(v_{\0}+v_{\1},Z(v_{\0}+v_{\1})))
\in
S(V_\0^*)\tens\PC.
\end{equation}
Now Fubini's formula gives:

\begin{equation}
\label{Decomp_Spf(X)}
\int_Vd_V(v)\exp(-\frac{1}{2}B(v,Xv))
=\int_{V_\0}d_{V_{\0}}(v_\0)\exp(-\frac{1}{2}B(v_\0,X_0v_\0))P(X-X_{0},v_\0).
\end{equation}
Since $X_0\in\UC^+$, $B(v_\0,X_0v_\0)$ is a positive definite quadratic form.
The integral on the right hand side is a Gaussian integral on $V_\0$. Thus 
$\Spf_\PC$ is an analytic function on
$\g{spo}(V)_\PC(\VC^+)$.

\subsubsection{$\Spf $ is a well defined analytic function on $\VC^+$}
\label{welldef}

 Let $\QC=\Lambda(\g{spo}(V)_{\1}^*)$. Let
$h:\g{spo}(V)_{\0}\hookrightarrow\big(\g{spo}(V)_{\0}\big)_{\QC}$ be the
canonical  embedding defined by $h(v)=v\tens 1$. Let $\Xi\in
 \big(\g{spo}(V)_{\1}\big)_{\QC}$ be the generic point of $\g{spo}(V)_{\1}$. We
put for $X\in\VC^+$:

\begin{equation}
\begin{split}
\phi(X)&=\Spf_{\QC}(h(X)+\Xi)\in \QC=\Lambda(\g{spo}(V)_{\1}^*)\\
&=\int_{V_\0}d_{V_{\0}}(v_\0)\exp(-\frac{1}{2}B(v_\0,X_{0}v_\0))P(X_1+\Xi,v_\0).
\end{split}
\end{equation}
 ($X=X_{0}+X_{1}$ with $X_{0}\in \UC^+$ and $X_{1}\in\g{so}(V_{\1})$; $P$ is
defined by (\ref{P}).) It defines a function 
\begin{equation}
\phi\in\CC^{\w}_{\g{spo}(V)}(\VC^+),
\end{equation}
such that for any near superalgebra $\PC$, any $X=Y+Z\in\g{spo}(V)_{\PC}(\VC^+)$ with 
$Y\in(\g{spo}(V)_{\0})_{\PC}(\VC^+)$ and $Z\in(\g{spo}(V)_{\1})_{\PC}$, 
$\phi(X)=\phi(Y)(Z)$.
\medskip

Since $\phi$ is defined by a Gaussian integral on $V_{\0}$ all its derivatives
along $\g{spo}(V)_{\0}$ are determined by derivation ``under the integral''.
Moreover the above integral is $\Lambda(\g{spo}(V)_{\1}^*)$-linear, hence for
any near superalgebra $\PC$  and $X\in\g{spo}(V)_{\PC}$,
$\phi(X)=\Spf_{\PC}(X)$.

\medskip

Now we put $\Spf=\phi\in\CC^{\w}_{\g{spo}(V)}(\VC^+)$ and call it the {\em
superPfaffian}. The preceding remark shows that for any near superalgebra $\PC$
and $X\in \g{spo}(V)_{\PC}(\VC^+)$ the expression $\Spf(X)$ is not ambiguous:
the value of the function $\Spf$ at $X$ is given by formula (\ref{eq:Spf}).

\medskip

We will see below (cf. section \ref{square}) that $\Spf^2=\Ber^-$.

\subsection{Holomorphic extension in the appropriate subset}\label{holo} 

Formula (\ref{eq:Spf}) is meaningful for $X\in\g{spo}(V\tens
\mathbb{C})_\PC$ with $\bb(X)\in\VC^+\times
\ii\,\g{spo}(V)_{\0}$ and it defines an holomorphic function on $\VC^+\times
\ii\,\g{spo}(V)$.

Indeed, let $X\in\g{spo}(V\tens \mathbb{C})_\PC$ with $\bb(X)\in\VC^+\times
\ii\,\g{spo}(V)_{\0}$. As before, let  
$X_0\in\UC^+\times\ii\,\g{sp}(V_\0)$ and $X_1\in\g{so}(V_\1\tens\mathbb{C})$
such that $\bb(X)=X_0+X_1$. The calculations of  section \ref{Spf_P} can be
reproduced here. The right hand side of formula (\ref{Decomp_Spf(X)}) is still a
Gaussian integral and therefore defines a complex analytic function on
$\g{spo}(V\tens
\mathbb{C})_\PC\big(\VC^+\times \ii\,\g{spo}(V)_{\0}\big)$. 

The same arguments as in section \ref{welldef} show that  formula (\ref{eq:Spf})
 defines an  holomorphic extension of $\Spf$ on $\VC^+\times
\ii\,\g{spo}(V)_{\0}$ still denoted by $\Spf$.

\subsection{Invariance}\label{Inv} Let $\PC$ be a near superalgebra. Let
$X\in\g{gl}(V)_{\PC}$. We denote by $X^*\in\g{gl}(V)_{\PC}$ the adjoint of $X$
defined by:
\begin{equation}
\label{adjoint}
\forall v,w\in V_{\PC},B(Xv,w)=B(v,X^*w).
\end{equation}
We have:
\begin{equation}
(X^*)^*=X.
\end{equation}

\medskip

Let $v\in V$, we denote by  $B^{\#}(v)$ the element of $V^*$ such that for any
$w\in V$, $B^{\#}(v)(w)=B(v,w)$. This defines an isomorphism $B^{\#}:V\to V^*$.
Moreover for $X\in\g{gl}(V)$ non zero and homogenous we denote by $^{t}X$ the
endomorphism of $V^*$ such that for any $\phi\in V^*$  non zero and homogenous
and any $v\in V$, $^{t}X(\phi)(v)=(-1)^{p(X)p(f)}\phi(Xv)$. Then:
\begin{equation}
\label{eq:X*} X^*=(B^{\#})^{-1}\,^{t}XB^{\#}.
\end{equation}

\medskip

We denote by $GL(V)_{\PC}$ the group of invertible elements of
$\g{gl}(V)_{\PC}$. 
Since $GL(V)_\PC\subset\g{gl}(V)_\PC$, the definition of $X^*$ is meaningful for
$X=g\in GL(V)_\PC$. For $g\in GL(V)_{\PC}$ we have from (\ref{eq:X*}):
$\Ber(g^*)=\Ber(g)$ and $\Ber_{(1,0)}(g^{*})=\Ber_{(1,0)}(g)$.

\medskip
We put:
\begin{displaymath}
SpO(V)_{\PC}
=\big\{g\in GL(V)_{\PC}\,/\,g^*=g^{-1} \big\}.
\end{displaymath}

\medskip

From the multiplicative property of $\Ber$ we get for $g\in SpO(V)_{\PC}$:
\begin{equation}
\label{eq:BerSpO}
\Ber(g)=\Ber_{(1,0)}(g)=\det(\bb(g)|_{V_{\1}})=\pm 1
\end{equation}

\begin{prop} \label{prop:Inv} Let $\PC$ be a near superalgebra. Let 
$X\in\g{spo}(V\tens\mathbb{C})_\PC\big(\VC^+\times\ii\,\g{spo}(V)\big)$ and 
$g\in GL(V)_\PC$, then:
\begin{align}
\label{eq:Inv}
\Spf(g^*X g)&=\Ber^{-1}_{(1,0)}(g)\Spf(X).
\end{align} In particular, for $g\in SpO(V)_\PC$, we have;
\begin{align}
\label{eq:InvSpO}
\Spf(g^{-1}X g)&=\det(\bb(g)|_{V_{\1}})\Spf(X).
\end{align}
\end{prop}

\begin{proof} We have by the formula of change of coordinates (cf.
\cite{Ber87}):
\begin{equation}
\begin{split}
\Spf(g^*Xg) & = \int_Vd_V(v) \exp\big(-\frac12B(v,g^*X gv)\big)\\
 & = \int_Vd_V(v) \exp\big(-\frac12B(gv,X gv)\big)\\
 & =\int_Vd_V(v)\Ber^{-1}_{(1,0)}(g) \exp\big(-\frac12B( v,X v)\big)\\
 & = \Ber^{-1}_{(1,0)}(g)\Spf(X).
\end{split}
\end{equation}

Now, assume that $g\in SpO(V)_\PC$. Then, by definition, $g^*=g^{-1}$ and
formula (\ref{eq:InvSpO}) follows from (\ref{eq:BerSpO}).

\end{proof}

\subsection{Action of differential operators}\label{Harm}

Let $V$ be a symplectic finite dimensional supervector space. We assume that
$\dim(V_{\1})$ is even.
\medskip 

Let $X\in\g{spo}(V)$ be homogeneous. We denote by $\de_{X}$ the derivation of
$\CC^\oo_{\g{spo}(V)}\big(\g{spo}(V)\big)$ such that for any homogeneous
$\psi\in\g {spo}(V_{\0})^*$:
\begin{equation}
\label{eq:der}
\de_X\psi=(-1)^{p(X)p(\psi)}\psi(X).
\end{equation}

The mapping $X\mapsto\de_X$ extends to an isomorphism between to
$S\big(\g{spo}(V)\big)$ and the superalgebra of differential operators with
constant coefficients on $\g{spo}(V)$.

Let $D\in S(\g{spo}(V))$. Since $\check\mu$ is linear and even on $\g{spo}(V)$,
we have (cf. formula (\ref{eq:checkmu}) for definition of $\check\mu$) for
any $X\in\g{spo}_{\PC}(V_{\0})$ and $v\in V_{\PC}$ where $\PC$ is a near superagebra:

\begin{equation}
\Big(\de_{D}\exp(\mu)\Big)(X,v)=\check\mu(D)(v)\exp(\mu)(X,v).
\end{equation} 

Moreover we recall that since $\Spf$  is defined by a Gaussian integral all its
derivatives are determined by derivation ``under the integral''. Thus, if
$X\in\g{spo}(V)_{\PC}(\VC^{+})$:

\begin{equation}
\label{deriv}
\begin{split}
\de_{D}\int_V d_V(v) \exp\big(\mu(X,v)\big)&=
\int_V d_V(v)\de_{D}\exp(\mu(X,v))\\
&=\int_V d_V(v)\check\mu(D)(v)\exp(\mu(X,v)).
\end{split}
\end{equation} 

It follows:

\medskip

\begin{prop} Let $D\in\ker(\check\mu)\subset S(\g{spo}(V))$. Then:
\begin{equation}
\de_{D}\,\Spf=0.
\end{equation}

\end{prop}

\bigskip

Let $X,Y\in\g{spo}(V)$. We put
\begin{equation}
K(X,Y)=\str(XY).
\end{equation}
It defines a non degenerate symmetric bilinear form on  $\g{spo}(V)$.

Let $(X_i)_{i\in I}$ be a basis of $\g{spo}(V)$ and $(X'_i)_{i\in I}$ the basis
of $\g{spo}(V)$ such that $K(X_i,X'_j)=\d_i^j$ ($\d_i^j$ is the Dirac symbol).
We put:

\begin{equation}
\square_K=\som_{i\in I} \de_{X'_i}\de_{X_i} \in S^2(\g {spo}(V)).
\end{equation}
It is an homogeneous differential operator of degree $2$ on $\g{spo}(V)$.

As a corollary of the  above proposition we obtain:

\begin{coro} Let $D\in \sdir_{k\in\mathbb N^*}S^k(\g{spo}(V))^{\g{spo}(V)}$. Thus
$\de_{D}$ is an $\g{spo}(V)$-invariant differential  operator on $\g{spo}(V)$
with constant coefficients and zero scalar term.

If $dim(V_\0)>0$:
\begin{equation}
\de_{D}\,\Spf=0.
\end{equation}

If $V=V_\1$:
\begin{equation}
deg(D)\not=\frac{dim(V)}{2}
\text{ or }
 D\in\Big(S^{\frac{dim(V)}{2} }\big(\g{so}(V)\big)\Big)^{O(V)}
\Rightarrow \de_{D}\,\Spf=0 .
\end{equation}

In particular, in all cases:
\begin{equation}
\square_K\,\Spf=0.
\end{equation}
\end{coro}

\begin{proof}[Proof] We will use the following lemma:
\begin{lemme}\label{LemInv}

For $k\geqslant 1$, if $dim(V_\0)>0$ or $V=V_\1$ and $k\not=dim (V)$:
\begin{equation}
\label{eq:Inv1} S^k(V^*)^{\g{spo}(V)}=\{0\}.
\end{equation}
If $V=V_1$ and $k=dim(V)$, we have:
\begin{equation}
\label{Inv2} S^k(V^*)^{\g{so}(V)}={\L}^k(V^*).
\end{equation}
\end{lemme}

\begin{proof}[Proof] Assume that   $dim(V_\0)>0$. Let $\PC$ be a near
superalgebra. Let $SpSO(V)_{\PC}$ be the connected component of $SpO(V)_{\PC}$.
Since for $v\in V_\0\setminus\{0\},\,SpSO(V)_{\PC}v=V_{\PC}$, the invariant
polynomials are constants and equality (\ref{eq:Inv1}) follows.

In case $V=V_\1$ cf. \cite{Wey46}.

\end{proof}

\medskip

Let $D\in S^k\big(\g{spo}(V)\big)^{\g{spo}(V)}$ $(k> 0)$, then $\check\mu(D)\in
S^{2k}(V^*)^{\g{spo}(V)}$.

Assume that $dim(V_\0)>0$ or $V=V_\1$ and $k\not=\frac{dim(V)}{2}$, then by
Lemma \ref{LemInv} we have $\check\mu(D)=0$.

Now assume that $V=V_\1$ and $k=\frac{dim(V)}{2}$. Then
$\check\mu(D)\in\L^{dim(V)} (V^*)$. In this case, if $\check\mu(D)$ is
$O(V)$ invariant, $\check\mu(D)=0$.

Now, the corollary follows from the proposition.

\end{proof}

\subsection{Taylor formula}\label{Taylor}
We still assume that $V$ is a symplectic finite dimensional supervector space
with $n=\dim(V_{\1})$ even.

\subsubsection{General case} 
Let $(P_{k})_{k\in\mathbb N}$ be an homogeneous (for parity) basis of
$S(V^*)$. We assume that $P_{k}$ is also homogeneous in degree as a polynomial.

 We define $c_{k}(X)$ as the coefficient of
$P_{k}$
in the expansion of $\exp$:
\begin{equation}\label{Taylorexp0}
\exp(\mu(X,v))=\som_{k\in\mathbb N}P_{k}(v)c_{k}(X),
\end{equation}
$c_{k}\in S(\g{spo}(V)^*)$. Then we define
 $\tilde {c}_{k}\in\CC^{\w}_{\g{spo}(V)}(\VC^+)$
by
\begin{equation}
\tilde{c}_{k}(X)=\int_{V}d_{V}(v) P_{k}(v)\exp\big(\mu(X,v)\big).
\end{equation}
We recall that $\Xi:\sdir_{k\in\mathbb N}S^{2k}(V^*)\to
\sdir_{k\in\mathbb N} S^{k}(\g{spo}(V))$ was defined in \ref{Xi}.

We put 
\begin{equation}
\de_{k}=\de_{\Xi(P_{k})}.
\end{equation}
We have:
\begin{lemme}\label{lem:c-de} For any near superalgebra $\PC$ and any
$X\in\g{spo}(V)_{\PC}(\VC^+)$:
\begin{equation}
\tilde c_{k}(X)=\big(\de_{k}\Spf\big)(X).
\end{equation}

\end{lemme}

\begin{proof}
We use formula (\ref{deriv}) to show:
\begin{equation}
\begin{split}
\de_{k}\int_{V}d_{V}(v)\exp(\mu(X,v))
&=\int_{V}d_{V}(v)\check\mu\big(\Xi(P_{k})\big)(v)\exp(\mu(X,v))\\
&=\int_{V}d_{V}(v)P_{k}(v)\exp(\mu(X,v))\\
&=\tilde c_{k}(X).
\end{split}
\end{equation}

\end{proof}

\medskip

\begin{prop} Let $\PC$ be a near superalgebra. Let $X,Y\in\g{spo}(V)_{\PC}$ 
such that 
 $\bb(Y)\in\g{so}(V_{1})$, and $\bb(X)\in\UC^{+}$. In this case $\bb(X+Y)\in\UC^+$. 
Taylor's formula for $\Spf$ reads:
\begin{equation}
\label{TaylorSpf}
\Spf(X+Y)=\som_{k\in\mathbb N}(-1)^{p(P_{k})}
c_{k}(Y)\tilde c_{k}(X).
\end{equation}
The sum converges as an analytic function in $X$.
\end{prop}
\begin{proof} We have 
\begin{displaymath}
\exp(\mu(X+Y,v))=\exp(\mu(Y,v))\exp(\mu(X,v)).
\end{displaymath} Then, we expand $\exp(\mu(Y,v))$ by formula
(\ref{Taylorexp0}):
\begin{displaymath}
\exp(\mu(Y,v))=\som_{k\in\mathbb N}P_{k}(v)c_{k}(Y)
=\som_{k\in\mathbb N}(-1)^{p(P_{k})}c_{k}(Y)P_{k}(v);
\end{displaymath}
 and integrate against $d_{V}(v)$ (since $\dim(V_{\1})$ is even, this operation is even 
 and thus commute 
 with multiplication by $c_{k}(Y)$ on the left).

\end{proof}

\bigskip

\subsubsection{Case $V=V_{\1}$.}  In the particular case where $V=V_{\1}$ and
$\Spf$ is the ordinary  pfaffian $\Pfaff$, we have the following
simplification. Since here the situation is purely algebraic, we can work with
$\mathbb C$ as ground field. 
We fix an oriented orthonormal basis $(f_{1},\dots,f_{n})$ of $V$. Let
$(\xi^{1},\dots,\xi^{n})$ be its dual basis. Then $(\xi^{J})_{J\in\{0,1\}^n}$
is a basis of $S(V^*)$. We define as above $c_{J},\, \tilde c_{J}$ and
$\de_{J}$.

For $J\in\{0,1\}^{n}$, we put
\begin{equation}
V_{J}=
\mathbb C j_{1}f_{1} + \dots + \mathbb C j_{n}f_{n}.
\end{equation}

For $J=(j_{1},\dots,j_{n})\in\{0,1\}^n$ we denote by
$J'=(j'_{1},\dots,j'_{n})\in\{0,1\}^n$ its complementary:
$j_{i}+j'_{i}=1$. We have
$V=V_{J}\oplus V_{J'}$. We denote by $p_{J}:V\to V_{J}$ the projection of $V$
onto $V_{J}$ with $\ker(p_{J})=V_{J'}$.

\medskip

For $|J|$ odd and $Y\in\g{so}(V)_{\PC}$, we have $c_{J}(Y)=0$. We now  consider the
case
$|J|$ even.

Since $(f_{1},\dots,f_{n})$ is an orthonormal oriented basis of $V$ the non-degenerate
symmetric bilinear form on $V_{\1}$ restricts to a non-degenerate symmetric
bilinear form on $V_{J}$. Let $(\xi^{1},\dots,\xi^{n})$ be the dual basis of 
$(f_{1},\dots,f_{n})$.  We give to $V_{J}$ the orientation defined by
$\xi^{J}.$ Let $1\le j_{1}<\dots<j_{r}\le n$ such that 
$\xi^{J}=\xi^{j_{1}}\dots\xi^{j_{r}}$. Then $(f_{j_{1}},\dots,f_{j_{r}})$ is an  
orthonormal oriented basis of $V_{J}.$

Let $Y\in\g{so}(V_{\1})$. We put:
\begin{equation}
\begin{split} Y_{J}:V_{J}&\rightarrow V_{J}\\
v&\mapsto Y_{J}(v)=p_{J}\big(Y(v)\big).
\end{split}
\end{equation}
We have $Y_{J}\in\g{so}(V_{J})$. The matrix of $Y_{J}$ in the basis
$(f_{j_{1}},\dots,f_{j_{r}})$ is obtained from the matrix of $Y$ in the basis
$(f_{1},\dots,f_{n})$ as the submatrix corresponding of rows and columns
$(j_{1},\dots,j_{r})$. 

We have:
\begin{equation}\label{cJPfaff1}
c_{J}(Y)=(-1)^{\frac{|J|(|J|-1)}2}\Pfaff(Y_{J}).
\end{equation}
For $J=(1,\dots,1)$ it is the definition of $\Pfaff$, and in the other cases
it follows (see \cite{MQ86}) by evaluating  $\exp(\mu(X,v))$ at
$\xi^{j'_{1}}=\dots=\xi^{j'_{n-r}}=0$.

We define  $\e(J,J')\in\{-1,1\}$ by the formula:
\begin{equation}
\e(J,J')\xi^{J}\xi^{J'}=\xi^{1}\dots\xi^{n};
\end{equation}
it is the signature of the permutation $(1,\dots,n)\mapsto(j_{1},\dots, j_{r},
j'_{1},\dots, j'_{n-r})$.

We obtain for $Y\in\g{so}(V_{\1})$:
\begin{equation}\label{cJPfaff2}
\begin{split}
c_{J'}(Y)
&=(-1)^{\frac{n(n-1)}2}\int_{V}d_{V}(v)(\xi^{1}\dots\xi^{n})(v)c_{J'}(Y)\\
&=(-1)^{\frac{n(n-
1)}2}\int_{V}d_{V}(v)\big(\e(J,J')\xi^{J}\xi^{J'}\big)(v)c_{J'}(Y)\\
&=(-1)^{\frac{n(n-1)}2}\e(J,J')\int_{V}d_{V}(v)\xi^{J}(v)\exp(\mu(Y,v))\\ 
&=(-1)^{\frac{n(n-1)}2}\e(J,J')\tilde c_{J}(Y).
\end{split}
\end{equation}
Since for $|J|$ even: 
$(-1)^{\frac{n(n-1)}2+\frac{|J|(|J|-1)}2+\frac{|J'|(|J'|-1)}2}=1=(-1)^{|J|}$,
 formula (\ref{TaylorSpf}) reads (cf. \cite{MQ86,Ste90}):
\begin{equation}
\label{TaylorPfaff} 
\Pfaff(X+Y)=\som_{J\in\{0,1\}^n/|J|\text{even}}\e(J,J')\Pfaff(X_{J})\Pfaff(Y_{J'}).
\end{equation}

\bigskip

\subsubsection{Case $X\in(\g{spo}(V)_{\0})_{\PC}$ and
$Y\in(\g{spo}(V)_{\1})_{\PC}$}\label{sec:Taylorex2}

We fix a symplectic oriented basis $(e_{i},f_{j})$
of $V$. Let $(x^{i},\xi^{j})$ its dual basis. Let 
$I=(i_1,\dots,i_{m})\in\mathbb N^m$ be a multiindice. We put:
\begin{equation}
x^{I}=(x^{1})^{i_{1}}\dots(x^{m})^{i_{m}}.
\end{equation}
Then:
\begin{displaymath}
\big(\xi^{J}x^{I}\big)_{(I,J)\in\mathbb N^m\times\{0,1\}^n},
\end{displaymath} is a basis  of $S(V^*)$.

For $I=(i_{1},\dots,i_{m})\in\mathbb N^m$ we put $|I|=i_{1}+\dots+i_{m}$ and   
 $I!=i_{1}!\dots i_{m}!$. Moreover we put:
\begin{eqnarray}
\frac{\de^{|I|}}{\de x^{I}}&=\frac{\de^{|I|}}{(\de x^{1})^{i_{1}}\dots(\de
x^{m})^{i_{m}}}\\
\frac{\de^{|J|}}{\de x^{J}}&=
\frac{\de^{|J|}}{(\de \xi^1)^{j_{1}}\dots(\de \xi^n)^{j_{n}}}.
\end{eqnarray}

\medskip

Let $\PC$ be a near superalgebra. Let $v\in V\tens V^*$ be the generic point of
$V$. Let $X\in\g{spo}(V)_{\PC}$.
We define $c_{I,J}(X)$ as the coefficient of
$\xi^J x^{I}$
in the Taylor formula:
\begin{equation}
\label{Taylorexp2}
\exp(\mu(X,v))=\som_{(I,J)\in\mathbb N^{m}\times\{0,1\}^{n}}{\xi^{J}x^{I}}c_{I,J}(X).
\end{equation}
(In particular, for $X\in\g{spo}(V)_{\PC}$ such that $\bb(X)\in\g{so}(V_{\1})$,
since $\mu(X,v)$ is nilpotent, the sum is finite.) It defines  $c_{I,J}\in
S^{\frac{|I|+|J|}2}(\g{spo}(V)^*)$.
We put:

\begin{equation}
\tilde c_{I,J}(X)=\int_{V}d_{V}(v)(\xi^{J}x^{I})(v)\exp\big(\mu(X,v)\big).
\end{equation}

and for $A\in\g{sp}(V_{\0})_{\PC}$:
 
\begin{equation}
\tilde c_{I}(A)=\int_{V_{\0}}d_{V_{\0}}(v_{\0})(x^{I})(v_{\0})
\exp\big(\frac12 B(v_{\0},Av_{\0})\big).
\end{equation}

\bigskip

To avoid confusing notations, in the rest of this paragraph we denote by $\bB$ the symplectic form on $V$.
\medskip
We put:
\begin{equation}
X=
\begin{pmatrix} A&0\\
0&D
\end{pmatrix}\text{ and }Y=
\begin{pmatrix} 0&B\\
C&0
\end{pmatrix}.
\end{equation}
with $A\in\g{sp}(V_{\0})_{\PC},$ $D\in\g{so}(V_{\1})_{\PC}$, $B\in
Hom(V_{\1},V_{\0})_{\PC}$ and $C=-B^*\in Hom(V_{\0},V_{\1})_{\PC}$. 
Where $B^*$ is defined by: 
$$
\forall v\in V_{\0}\tens \PC_{\0},\ 
\forall w\in V_{\1}\tens \PC_{\1},\ 
 \bB(B^*v,w)=\bB(v,Bw)
$$

\medskip

Let $v_{\0}$ be the generic point of $V_{\0}$, $v_{\1}$ be the generic point of 
$V_{\1}$ and $v=v_{\0}+v_{\1}$ be the generic point of $V$. 
 We have:
  \begin{equation}
\mu(X,v)=-\frac12\bB(v_{\0},Av_{\0})-\frac12\bB(v_{\1},Dv_{\1}).
\end{equation}

  Assume that $\bb(X)\in\VC^+$ that means $\bb(A)\in\UC^+$.  Then, with the notations
of the preceding subsection, for $|J|$ even, we have $(-1)^{\frac{|J'|(|J'|-
1)}2+\frac{n(n-1)}2}=(-1)^{\frac{|J|(|J|-1)}2}$, and so:
\begin{equation}
\begin{split}
\tilde c_{I,J}(X)&=
\int_{V_{\0}}d_{V_{\0}}(v_{\0})x^{I}(v_{\0})\exp(-\frac12\bB(v_{\0},Av_{\0}))
\int_{V_{\1}}d_{V_{\1}}(v_{\1})
\xi^{J}(v_{\1})\exp(-\frac12\bB(v_{\1},Dv_{\1}))\\
&=(-1)^{\frac{|J|(|J|-1)}2}\e(J,J')
\tilde c_{I}(A)\Pfaff(D_{J'}).
\end{split}
\end{equation}

\bigskip Now, we explicit $c_{I,J}(Y)$. Let us introduce some notations.

Since $B=-C^*$ we
have: 
\begin{displaymath}
\begin{split}
\mu(Y,v)&=-\frac12\big(\bB(v_{\0},Bv_{\1})+\bB(v_{\1},Cv_{\0})\big)\\
&=-\frac12\big(\bB(B^*v_{\0},v_{\1})+\bB(v_{\1},Cv_{\0})\big)\\
&=-\frac12\big(-\bB(Cv_{\0},v_{\1})+\bB(v_{\1},Cv_{\0})\big)\\
&=-\bB(v_{\1},Cv_{\0}).
\end{split}
\end{displaymath}

\medskip

Let $(I,J)\in\mathbb N^m\times\{0,1\}^{n}$. We denote by $C_{J,I}$ the $|J|\times
|I|$ matrix obtained from $C$ by keeping $j_{k}$ times the $k$-th line of
$C$ (in other words we keep the
lines $(j_{1},\dots,j_{r})$) and $i_{k}$ times the $k$-th column of $C$. 

\emph{Example:} Assume that $(m,n)=(3,4)$. Let:
\begin{displaymath} C=
\begin{pmatrix}
\a_{1}&\a_{2}&\a_{3}\\
\be_{1}&\be_{2}&\be_{3}\\
\gamma_{1}&\gamma_{2}&\gamma_{3}\\
\d_{1}&\d_{2}&\d_{3}
\end{pmatrix}.
\end{displaymath} Let $I=(2,0,1)$ and $J=(0,1,1,1)$. Then:
\begin{displaymath} C_{J,I}=
\begin{pmatrix}
\be_{1}&\be_{1}&\be_{3}\\
\gamma_{1}&\gamma_{1}&\gamma_{3}\\
\d_{1}&\d_{1}&\d_{3}\\
\end{pmatrix}.
\end{displaymath}

\medskip Let $r\in\mathbb N$. We denote by $\mathfrak{S}_{r}$ the group of
permutations of $\{1,\dots,r\}$. We denote by $\phi_{r}$ the $r$-multilinear form on
the $r\times r$ matrix antisymmetric in the lines, symmetric in the columns
defined for  $M=(a_{i,j})_{1\le i,j\le r}$ with $a_{i,j}\in\PC_{\1}$ by:
\begin{equation}
\phi_{r}(M)=\som_{\s\in\g{S}_{r}} a_{1,\s(1)}\dots a_{r,\s(r)}
\end{equation}

 \medskip
 
 In the sequel we put for $C\in Hom(V_{\0},V_{\1})\tens\PC_{\1}$ and
$I,J\in\mathbb N^m\times \{0,1\}^{n}$:
 \begin{equation}
c_{I,J}(C)=
\begin{cases} 0&\text{ if }|I|\not=|J|;\\
(-1)^{\frac{|J|(|J|-1)}2}\phi_{|J|}(C_{J,I})&\text{ if }|I|=|J|.
\end{cases}
\end{equation}

\medskip

With this notations we have:

\begin{equation}\label{cIJ}
c_{I,J}(Y)=c_{I,J}(C).
\end{equation}
 
\medskip

\emph{Example:} We take the preceding example with $(m,n)=(3,4)$,
$Y=\begin{pmatrix} 0&B\\
C&0
\end{pmatrix}$ with $C$ as above,  $I=(2,0,1)$ and $J=(0,1,1,1)$. Then:
\begin{displaymath}
c_{I,J}(Y)=c_{I,J}(C)=-\phi_{3}(C_{J,I})=-2(\be_{1}\gamma_{1}\d_{3}
+\be_{1}\gamma_{3}\d_{1}+\be_{3}\gamma_{1}\d_{1})
\end{displaymath}

\bigskip

Now formula (\ref{TaylorSpf}) gives:

\begin{prop} Let $\PC$ be any near superalgebra. Let $
\begin{pmatrix} A&B\\
C&D
\end{pmatrix}\in\g{spo}(V)_{\PC}(\VC^+)$, then:
\begin{equation}
\label{TaylorSpf2}
\Spf
\begin{pmatrix} A&B\\
C&D
\end{pmatrix}=\som_{(I,J)\in\mathbb N^m\times\JC_{n}\,/\, |I|=|J|\text{ even }}
(-1)^{\frac{|J|(|J|-1)}2}\e(J,J')c_{I,J}(C)\Pfaff(D_{J'})
\tilde c_{I}(A).
\end{equation}
\end{prop}

\bigskip

\subsection{$\Spf(-X^{-1})$ and $\Spf^2$}

\subsubsection{Some formulas}

Let $W$ be a supervector space. Let $\bw\in (V\tens W^*)_{\0}$. Let $\PC$ be any
near superalgebra. For $w\in W_{\PC}$, $\bw(w)$ belongs to $ V_{\PC}$ and 
$v\mapsto
 B(v,\bw(w))$ is linear on $V$ while $v\mapsto B(v,Xv)$ is quadratic.
Thus, as in proof of theorem \ref{th:Spf} (cf. section \ref{correct}),  for any
near superalgebra $\PC$ we define an analytic function $\psi_{\PC}$ on
$\g{spo}(V)_{\PC}(\VC+)\times W_{\PC}\times\mathbb C$ by the formula
($(X,w,\l)\in\g{spo}(V)_{\PC}(\VC+)\times W_{\PC}\times\mathbb C$):
\begin{equation}
\psi_{\PC}(X,w,\l)=\int_{V}d_{V}(v)\exp\big(\mu(X,v)+\l B(v,\bw(w)\big).
\end{equation}

Then, as in the proof of theorem \ref{th:Spf}, we can prove that there is a
function  $\psi$ on $\g{spo}(V)\times W\times \mathbb C$ defined on $\VC^+\times
W_{\0}\times\mathbb C$ such that for any near superalgebra $\PC$ and 
$\forall(X,w,\l)\in\g{spo}(V)_{\PC}(\VC^+)\times W_{\PC}\times\mathbb C$,
$\psi(X,w,\l)=\psi_{\PC}(X,w,\l)$.

Now, we put:
\begin{displaymath}
\Spf^{\bw}_{\l}(X,w)=\psi(X,w,\l).
\end{displaymath}
Thus, for $\l\in\mathbb C$, $\Spf^{\bw}_{\l}\in\CC^{\oo}_{\g{spo}(V)\times
W}(\VC^+\times
W_{\0}).$

\medskip

\begin{lemme}(cf. \cite{MQ86} for the case $V=V_{\1}$.) Let
$X\in\g{spo}(V)_{\PC}(\VC^+)$ be invertible (since $\bb(X)|_{V_{\0}}$ is
already invertible, it means $\bb(X)|_{V_{\1}}$ invertible)
and $w\in W_{\PC}$. We have:
\begin{equation}
\label{eq:Spf-bw}
\Spf^{\bw}_{\l}(X,w)=\Spf(X)\,\exp(\frac{-\l^2}2B(\bw(w),X^{-1}\bw(w)))
\end{equation}

\end{lemme}

\begin{proof} 
Since $X^*=-X$, we have:
\begin{equation}
\begin{split}
\mu(X,v)+\l B(v,\bw(w))&=-\frac12 B(v,Xv)+\l B(v,\bw(w))\\
&=-\frac12 B(v-\l X^{-1}\bw(w),X(v-\l X^{-1}\bw(w)))\\
&\hskip 1cm
+\frac{\l^2}2 B(X^{-1}\bw(w),\bw(w)).
\end{split}
\end{equation}

 We put:
\begin{equation}
\phi(X,\l,w)=\int_{V}d_{V}(v)\exp(-\frac12 B(v-\l X^{-1}\bw(w),X(v-\l X^{-1}\bw(w)))).
\end{equation}
It is an analytic function on $\VC^+\times \mathbb C\times V$. Since $d_{V}(v)$ is 
invariant by translations, we have on $\VC^+\times \mathbb R\times V$: 
$\phi(X,\l,w)=\Spf(X)$.
By uniqueness of analytic continuation, it follows that for any $\l\in\mathbb
C$, $\phi(X,\l,w)=\Spf(X)$. 
\end{proof}

\bigskip

Applying this lemma for various particular values of  $W$ and $\bw$ we will
obtain some useful formulas. First, we take $W$ to be a supervector space
isomorphic to $V$ and $\bw\in V\tens W^*$ be an isomorphism $W\to V$. 

\medskip

Let $D\in S(V)$. It defines a differential operator $\de_{D}$ on $V$.
We put:
\begin{equation}
\overline {c}_{D}(X)=\Big(\de_{D}\exp(\mu)\Big)(X,0).
\end{equation}
and for $P\in S(V^*)$ and $X\in\g{spo}(V)_{\PC}(\VC^+)$:
\begin{equation}
\tilde {c}_{P}(X)=\int_{V}d_{V}(v)P(v)\exp(\mu(X,v)).
\end{equation}

 We recall (cf. section \ref{Inv}) that $B^{\#}:V\to V^*$
 is an  isomorphism. It extends to an isomorphism of algebras $B^{\#}:S(V)\to S(V^*)$. 
 We put for $D\in S(V)$:
 \begin{equation}
D^{\#}=B^{\#}(D).
\end{equation}

\medskip

Let $v\in V\tens V^*$ be the generic point of $V$. For $D\in\ S^k(V)$ we have:
\begin{equation}
\de_{D}v^{k}=k!D;
\end{equation}
and:
\begin{equation}
\de_{\bw^{-1}(D)}\big(B(v,\bw)^k\big)=(-1)^k k!D^{\#}(v).
\end{equation}

\medskip

{\em Example:} Let $(e_{1},\dots,e_{m},f_{1},\dots,f_{n})$ be an homogeneous basis of $V$. 
Let $(x^{1},\dots,x^{m},\xi^{1},$ $\dots,\xi^{n})$ be its dual basis. 
Let $(I,J)\in\mathbb N^m\times\{0,1\}^n$. We put 
$D_{I,J}=e_{1}^{i_{1}}\dots e_{m}^{i_{m}}f_{1}^{j_{1}}\dots f_{n}^{j_{n}}$. 
We have for $X\in\g{spo}(V)_{\PC}$ ($\de_{e_{i}}=\frac{\de}{\de x^{i}}$ and 
$\de_{f_{j}}=-\frac{\de}{\de \xi^{j}}$):
\begin{align*}
\de_{D_{I,J}}\xi^{J}x^{I}
 &=(-1)^{|J|}
 \big(\frac{\de}{\de x^1}\big)^{i_{1}}\dots
\big(\frac{\de}{\de
x^{m}}\big)^{i_{m}}\big(\frac{\de}{\de \xi^1}\big)^{j_{1}}\dots
\big(\frac{\de}{\de \xi^{n}}\big)^{j_{n}}
\xi^{J}x^{I}\\
&={I!}(-1)^{\frac{|J|(|J|+1)}2}
\end{align*}

Then:

\begin{align*}
\overline {c}_{D_{I,J}}(X)&=I! (-1)^{\frac{|J|(|J|+1)}2}c_{I,J}(X)
\\
\text{and}\quad
\exp(\mu(X,v))&=\som_{(I,J)\in\mathbb N^m\times\{0,1\}^{n}}
\frac{\xi^{J}x^{I}}{I!(-1)^{\frac{|J|(|J|+1)}2}}\overline c_{D_{I,J}}(X).
\end{align*}

\medskip

We obtain:

\begin{coro} For any $X\in
\g{spo}(V)_{\PC}(\VC^+)$ invertible, and $D\in S^{2k}(V)$ (if $D\in S^{2k+1}(V)$
$\tilde c_{D^{\#}}(X) =c_{D}(-X^{-1})=0$):
\begin{equation}\label{eq:deSpf}
\tilde c_{D^{\#}}(X) 
=\overline c_{D}(X^{-1})\Spf(X).
\end{equation}
and for $\Re (\l^2)<0$ (in this case $\l^2X^{-
1}\in\g{spo}(V)_{\PC}(\VC^+)\times\ii\g{spo}(V)_{\PC}$):
\begin{equation}
(-\l^2)^{\frac{m-n-2k}2}(-1)^k\overline c_{D}(X)\overline c_{D^{\#}}(X)
=\tilde c_{D^\#}(\l^2 X^{-1})\,\Spf(X).
\end{equation}
\end{coro}
\medskip

{\em Remark:} We point out that in case $V=V_{\0}$ formula (\ref{eq:deSpf}) is Wick 
formula (cf. for example \cite{GJ81}).

\medskip
\begin{proof}   Let $v\in V\tens V^*$ be the generic point of $V$. Let $D\in S^{2k}(V)$.
Then:

\begin{displaymath}
\de_{\bw^{-1}(D)}\exp(\l B(v,\bw))=(-\l)^{2k}D^{\#}(v)\exp(\l B(v,\bw)).
\end{displaymath} 

\medskip

Now, we apply $ \frac 1{\l^{2k}}\de_{\bw^{-1}(D)}$
to  equality (\ref{eq:Spf-bw}) and then, we take the value at
$(X,0)$. Since,  $\overline c_{D}(\l^2X^{-
1})=\l^{2k}\overline c_{D}(X^{-1})$, formula (\ref{eq:deSpf}) follows.

 \bigskip
 
For  the second formula,  we need an auxiliary result.

Consider the application:
\begin{displaymath}
 \phi\mapsto\FC_{\l}(\phi)=\int_{W}d_{W}(w)\int_{V}d_{V}(v)\phi(v)\exp(\l
B(\bw(w),v)).
\end{displaymath}
It is defined for $\phi\in\CC^{^\oo}_{V}(V_{\0})$ such that for any $w\in W_{\PC}$,
 $v\mapsto\phi(v)\exp(\l
B(\bw(w),v))$ and all its derivatives is rapidly decreasing on $V$. It is linear and $\phi(0)=0$  implies that $\FC_{\l}(\phi)=0$. Thus, there is  
 $K_{\l}\in\mathbb C$ such that $\FC_{\l}(\phi)=K_{\l}\phi(0)$. To find $K_{\l}$
it is enough to consider a particular $\phi$. For example
$\phi(v)=\exp(\mu(X,v))$ ($\phi=\exp(\check\mu(X))$) for some
$X\in\g{spo}(V)_{\PC}(\VC^+)$ fixed. In this case formula (\ref{eq:Spf-bw}) shows that:
\begin{displaymath}
\FC_{\l}\big(\exp(\check\mu(X))\big)=\int_{V}d_{V}(v)\exp(\mu(X,v))
\int_{W}d_{W}(w)\exp(\frac{-\l^2}2B(\bw(w),X^{-1}\bw(w))).
\end{displaymath}
hypothesis $\Re(\l^2)<0$ implies that $\l^2 X^{-1}\in\g{spo}(V)_{\PC}(\VC^+)$ and thus,
 the above integral converges. More precisely, since $\exp(\mu(X,0))=1$, we have:
 \begin{equation}\label{eq:K}
K_{\l}=\FC_{\l}\big(\exp(\check\mu(X))\big)=\Spf(X)\Spf(\l^2X^{-1})
.
\end{equation}
Since it does not depends on $X$ it is enough to take 
\begin{displaymath}
X=\begin{pmatrix}
\vbox{\offinterlineskip\halign{\hglue 5pt#\hglue 5pt&#&
\hglue 5pt#\hglue 5pt\cr $
\begin{matrix}J_{2}&&0\\
&\ddots&\\
0&&J_{2}
\end{matrix}$ &\vrule height 18pt depth 12pt& $
\begin{matrix}0&&0\\
&\ddots&\\
0&&0
\end{matrix}$\cr
\noalign{\hrule} $
\begin{matrix}0&&0\\
&\ddots&\\
0&&0
\end{matrix}$\hglue 6pt&
\vrule height 18pt depth 13pt&\hglue 11pt$
\begin{matrix}J_{2}&&0\\
&\ddots&\\
0&&J_{2}
\end{matrix}$\cr}}
\end{pmatrix}
\end{displaymath}
where $J_{2}=\begin{pmatrix}0&-1\\
1&0
\end{pmatrix}$. We obtain:
\begin{equation}
K_{\l}=(-\l^{2})^{\frac{n-m}2}.
\end{equation}

\bigskip
Now, we multiply both sides of (\ref{eq:Spf-bw}) by $D^{\#}(\bw(w))$ and then  integrate on
$W$  against $d_{W}(w)$.

 For the left hand side we have
 (the first equality is obtained by integration by parts):
 
 \begin{displaymath}
\begin{split}
\int_{W}d_{W}(w)&D^{\#}(w)\int_{V}d_{V}(v)\exp(\mu(X,v)+\l B(v,\bw(w)))\\
&=
(-\l)^{-2k}\int_{W}d_{W}(w)\int_{V}d_{V}(v)\Big(\de_{D}\exp(\mu)\Big)(X,v)
\exp(\l B(v,\bw(w)))\\
&=(-\l)^{-2k}\FC_{\l}\Big(\de_{D}\exp(\check\mu(X))\Big)\\
&=(-\l)^{-2k}(-\l^2)^{\frac{n-m}2}\Big(\de_{D}\exp(\check\mu(X))\Big)(0)\\
&=(-\l^2)^{\frac{n-m-2k}2}(-1)^k\overline c_{D}(X).
\end{split}
\end{displaymath}
\end{proof}

\bigskip

As a corollary of equation (\ref{eq:Spf-bw}) we also have:

\begin{coro}\label{coro1} Let $\PC$ be a near superalgebra. 

Let
$A\in\g{sp}(V_{\0})_{\PC}(\UC^+)$, $B\in Hom(V_{\1},V_{\0})\tens\PC_{\1}$ and
$C\in Hom(V_{\0},V_{\1})\tens\PC_{\1}$ such that $\begin{pmatrix}
A&B\\
C&0
\end{pmatrix}
\in\g{spo}(V)_{\PC}.$ Then for any $J\in\{0,1\}^{n}$
even:
\begin{equation}
\label{eq1}
\begin{split}
\som_{I\in\mathbb N^m\,/|I|=|J|}(-i)^{|J|}
\tilde c_{I}(A)c_{I,J}(C)&=\frac1{\sqrt%
{\det(A)}}c_{J}(CA^{-1}B)\\
&=\frac{(-1)^{\frac{|J|(|J|-
1)}2}}{\sqrt{\det(A)}}\Pfaff((CA^{-1}B)_{J}).
\end{split}
\end{equation}

\end{coro}

\begin{proof} In this proof to avoid confusing notations, as in subsection
\ref{sec:Taylorex2}
we denote by $\bB$ the symplectic form on $V$.

 For $v\in (V_{\0})_{\PC}$ and $w\in V_{\1}$. we have:
 \begin{equation}
\bB(v,Bw)=-\bB(Cv,w)=\bB(w,Cv).
\end{equation}
Thus if $\PC=S(V^*)$, $v=\som_{i}e_{i}x^{i}$ (resp. $w=\som_{j}f_{j}\xi^{j}$) is
the generic point of $V_{\0}$ (resp. $V_{\1}$), the coefficient of
$\xi^{J}x^{I}$  in $\exp(-\bB(v,Bw))=\exp\big(-\frac12(\bB(v,Bw)+\bB(w,Cv))\big)$ is
$c_{I,J}(C)$ (cf. (\ref{cIJ})).
On the other hand, we have for $w\in (V_{\1})_{\PC}$:
\begin{equation}
\bB(Bw,A^{-1}Bw)=-\bB(w,CA^{-1}Bw);
\end{equation}

Now, look at equation (\ref{eq:Spf-bw}) with $V=V_{\0}$, $W=V_{\1}$, $X=A$,
and:
\begin{equation}
\bw(w)=Bw.
\end{equation}
 Then we apply $\big(\frac{\de}{\de
\xi^{n}}\big)^{j_{n}}\dots \big(\frac{\de}{\de \xi^{1}}\big)^{j_{1}}$ to
(\ref{eq:Spf-bw}) and taking value at $(X,0)$. 
Then equality (\ref{eq1}) follows by multiplying each side by $\frac{(-\ii)^{|J|}}{\l^{|J|}}$.
\end{proof}

\medskip

Similarly, we obtain:

\begin{coro}\label{coro2} Let $\PC$ be a near superalgebra. 

Let $D\in\g{so}(V_{\1})_{\PC}$ invertible, $B\in Hom(V_{\1},V_{\0})\tens\PC_{\1}$
and $C\in Hom(V_{\0},V_{\1})\tens\PC_{\1}$  such that $\begin{pmatrix}
0&B\\
C&D
\end{pmatrix}
\in\g{spo}(V)_{\PC}$. Then for any $I\in\mathbb N^m$:
\begin{equation}
\label{eq2}
\som_{J\in\JC_{n}\,/|I|=|J|}
\e(J,J')(-1)^{\frac{|J|(|J|-1)}2}\Pfaff(D_{J'})c_{I,J}(C) =\Pfaff(D)c_{I}(BD^{-1}C).
\end{equation}

\end{coro}

\begin{proof} Here we apply formula (\ref{eq:Spf-bw}) with $V=V_{\1}$,
$W=V_{\0}$,  $X=D$ and 
\begin{displaymath}
\bw(w)=Cw.
\end{displaymath}
Now, we  apply
\begin{displaymath}
 \frac 1{I!}\frac{\de^{i_{1}}}{\de y_{1}^{i_{1}}}\dots \frac{\de^{i_{m}}}{\de
y_{m}^{i_{m}}}
\end{displaymath}
 to (\ref{eq:Spf-bw}) and taking value at $(X,0)$.

Then, using (\ref{cJPfaff1}) and (\ref{cJPfaff1}), equality (\ref{eq1}) follows from 
multiplication of each side by $\frac{(-1)^{\frac{|J|}2}}{\l^{|J|}}$.
\end{proof}

\medskip

\subsubsection{Evaluation of $\Spf^2$}\label{square}

\begin{prop}

We have in $\CC^{\w}_{\g{spo}(V)}(\VC^+)$.
\begin{equation}
\label{eq:square}
\Spf^2=\Ber^-.
\end{equation}
It is equivalent to say that for any near superalgebra $\PC$ and any
$X\in\g{spo}(V)_\PC(\VC^+)\times\ii\,\g{spo}(V)_\PC$ we have:
\begin{equation}
\Spf^2(X)=\Ber^-(X).
\end{equation}
\end{prop}

\begin{proof} Proposition
\ref{prop:Inv} with $g=X^{-1}$  gives:
\begin{equation}
\Spf(-X^{-1})=\Ber^{-1}_{(1,0)}(X^{-1})\Spf(X).
\end{equation}

Now, since on $\VC^+$, $\Ber_{(1,0)}^{-1}=\Ber^{-}$, the result follows from multiplying  
both sides by $\Ber^-(X)\Spf(X)$ and using
formula (\ref{eq:K}) with $\l=\ii$.

\end{proof}

\subsection{Product formulas}

\subsubsection{Case $X\in\g{spo}(V)_{\PC}(\UC^+)$}

Let $X=
\begin{pmatrix}A&B\\
C&D
\end{pmatrix}$. We put $X^*=
\begin{pmatrix}A^*&C^*\\
B^*&D^*
\end{pmatrix}$. Then, $X\in\g{spo}(V)_{\PC}\Leftrightarrow X^*=-X$. Moreover,
$X\in\g{spo}(V)_{\PC}(\UC^+)$ implies that $A$ is invertible. It follows that:
\begin{equation}
\begin{split} A^*&=-A\\
D^*&=-D\\
C^*&=-B
\end{split}
\end{equation}
Thus:
\begin{equation}
\begin{split}
\begin{pmatrix}A&B\\
C&D
\end{pmatrix}&=
\begin{pmatrix}1&0\\
CA^{-1}&1
\end{pmatrix}
\begin{pmatrix}A&0\\
0&D-CA^{-1}B
\end{pmatrix}
\begin{pmatrix}1&A^{-1}B\\
0&1
\end{pmatrix},\\
&=
\begin{pmatrix}1&A^{-1}B\\
0&1
\end{pmatrix}^*
\begin{pmatrix}A&0\\
0&D-CA^{-1}B
\end{pmatrix}
\begin{pmatrix}1&A^{-1}B\\
0&1
\end{pmatrix}.
\end{split}
\end{equation}

Since $\Ber_{(1,0)}
\begin{pmatrix}1&A^{-1}B\\
0&1
\end{pmatrix}=1$, 
formula (\ref{eq:Inv}) imply  
\begin{equation}
\label{eq:SpfDecomp}
\Spf
\begin{pmatrix}A&B\\
C&D
\end{pmatrix}=
\Spf
\begin{pmatrix}A&0\\
0&D-CA^{-1}B
\end{pmatrix}=\frac{\Pfaff(D-CA^{-1}B)}{\sqrt{\det(A)}}.
\end{equation}

\medskip

We explain the relation with formula (\ref{TaylorSpf2}).

Taylor formula (\ref{TaylorPfaff}) for $\Pfaff(D-CA^{-1}B)$ gives:
\begin{equation}
\Pfaff(D-CA^{-1}B)=\som_{J\in\{0,1\}^{n}/|J|\text{ even}}
\e(J,J')\Pfaff((-CA^{-1}B)_{J})\Pfaff(D_{J'}).
\end{equation}
Compatibility with formula (\ref{TaylorSpf2}) follows from formula (\ref{eq1}) and that for
 $|J|$ even  $(-\ii)^{|J|}=(-1)^{\frac{|J|(|J|-1)}2}$.

\bigskip

\subsubsection{Case $X\in\g{spo}(V)_{\PC}(\UC^+)$ and invertible.} As before we
put $X=
\begin{pmatrix} A&B\\
C&D
\end{pmatrix}$ but now we assume moreover that $D$ is invertible. We have:

\begin{equation}
\begin{split}
\begin{pmatrix}A&B\\
C&D
\end{pmatrix}&=
\begin{pmatrix}1&BD^{-1}\\
0&1
\end{pmatrix}
\begin{pmatrix}A-BD^{-1}C&0\\
0&D
\end{pmatrix}
\begin{pmatrix}1&0\\
D^{-1}C&1
\end{pmatrix},\\
&=
\begin{pmatrix}1&0\\
D^{-1}C&1
\end{pmatrix}^*
\begin{pmatrix}A-BD^{-1}C&0\\
0&D
\end{pmatrix}
\begin{pmatrix}1&0\\
D^{-1}C&1
\end{pmatrix}.
\end{split}
\end{equation}

Since $\Ber_{(1,0)}
\begin{pmatrix}1&0\\
D^{-1}C&1
\end{pmatrix}=1$, 
formula (\ref{eq:Inv}) imply  
\begin{equation}
\Spf
\begin{pmatrix}A&B\\
C&D
\end{pmatrix}=
\Spf
\begin{pmatrix}A-BD^{-1}C&0\\
0&D
\end{pmatrix}=\frac{\Pfaff(D)}{\sqrt{\det(A-BD^{-1}C)}}.
\end{equation}

\medskip Now, compatibility with formula (\ref{TaylorSpf2}) comes from
formula (\ref{TaylorSpf}) for $\frac 1{\sqrt{det(A-BD^{-1}C)}}$ and equation (\ref{eq2}).

\medskip

\subsection{Homogeneity}\label{homog}

\begin{coro} The superPfaffian is an homogeneous function of degree $\frac{n-m}
2$ on $\VC^+\times\ii\,\g{spo}(V)$.
\end{coro}
\begin{proof} Let $\l>0$. We denote by $M_\l\in GL(V_{\0})\times GL(V_{\1})$ the
homothecy with ratio $\l$ on $V$. For any near superalgebra $\PC$, any
$X\in\g{spo}(V\tens \mathbb{C})_\PC(\VC^+\times\ii\,\g{spo}(V))$, we have:
\begin{equation}
\l X=M_{\sqrt\l} \, X\, M_{\sqrt\l}
\end{equation}
The corollary follows from proposition \ref{prop:Inv} and the equalities 
\begin{equation}
M_{\sqrt\l}^*=M_{\sqrt\l}
\ \text{ and }\  
\Ber_{(1,0)}(M_{\sqrt\l})=\big(M_{\sqrt\l}\big)=\l^{\frac{m-n}2}.
\end{equation}
\end{proof}

\section{SuperPfaffian II : a generalized  function}

In this section we define $\Spf$ as a generalized function on $\g{spo}(V)$ by
the formula:
\begin{equation}
\label{eq:SpfII}
\Spf(X)=\ii^{\frac{m-n}2}\int_Vd_V(v)\exp(-\frac\ii2B(v,Xv))
\end{equation}

The meaning of this formula is the following. For any smooth compactly supported
distribution  $t$ on $\g{spo}(V)$:
\begin{equation}
\label{def-gene}
\int_{\g{spo}(V)}t(X)\Spf(X)=\ii^{\frac{m-n}2}\int_Vd_V(v)\int_{\g{spo}(V)}t(X%
)\exp(-\frac\ii2B(v,Xv)).
\end{equation}
This means:
\begin{enumerate}
  \item we evaluate  the integral on $\g{spo}(V)$;
  \item the resulting function is rapidly decreasing on $V$ (cf. section
\ref{Spf:gene});
  \item we evaluate the integral on $V$.
\end{enumerate}

This generalized function on the supermanifold $\g{spo}(V)$ coincides on $\VC^+$
with the superPfaffian defined in the preceding section (cf. subsection
\ref{analytic} for a proof).

\subsection{A well defined generalized function}\label{Spf:gene} In this section
we prove  that formula (\ref{eq:SpfII}) defines a generalized function on
$\g{sp}(V_\0),$ with values in a finite dimensional subspace of
$S\big(\big(\g{so}(V_\1)\oplus\g{spo}(V)_\1\big)^*\big)$ in a sense analog to
formula (\ref{def-gene}). In particular this ensures  that formula
(\ref{eq:SpfII}) defines generalized function on $\g{spo}(V)$. The point is to
show that $\int_{\g{spo}(V)}t(X)\exp(-\frac\ii2B(v,Xv))$ is rapidly decreasing
on $V$.

\medskip

Let us precise some notations.

Let $v\in V\tens V^*$ be the generic point of $V$.  We denote by
$\widetilde\mu\in \g{spo}(V)^*\tens S^2(V^*)$ the polynomial of degree $2$ on
$V$ with values in  $\g{spo}(V)^*$ such that for any near superalgebra $\PC$ and
any $X\in\g{spo}(V)_\PC$:
\begin{equation}
\label{eq:moment}
\widetilde\mu(v)(X)=\mu(X,v).
\end{equation}
In particular for $u\in V_\PC$, $\widetilde\mu(u)\in\g{spo}(V)_{\PC}^*$.

\medskip

Let $\r$ be a smooth compactly supported distribution on $\g{sp}(V_\0)$. We
denote by $\what\r$ its Fourier transform. It is a smooth rapidly decreasing
function on $\g{sp}(V_\0)^*$ (in sense of Schwartz) which is defined for
$f\in\g{sp}(V_\0)^*$ by the formula:
\begin{equation}
\label{eq:Fourier}
\what\r(f)=\int_{\g{sp}(V_\0)}\r(X)\exp(-\ii f(X)),
\end{equation}

\medskip

We fix $J\in \UC^+$. Then $B(v,Jv)$ is a positive definite  quadratic form on
$V_\0$. We put for $u\in V_\0$, $\|u\|=\sqrt{\frac12 B(u,Ju)}$. 

We fix a norm  $N'$ on $\g{spo}(V)_\0$ and denote $N$ the associated norm on
$\g{spo}(V)_\0^*$. For $f\in\g{spo}(V)_\0^*$ we have:
\begin{equation}
N(f)=\sup_{Y\in\g{spo}(V)_\0\setminus\{0\}}\frac{|f(Y)|}{N'(Y)}.
\end{equation}
Thus we have for $u\in V_\0$:
\begin{equation}
\label{eq:MinN} N(\widetilde\mu(u)|_{\g{sp}(V_{\0})})\geqslant 
\frac{|\widetilde\mu(u)(J)|}{N'(J)}= \frac{\|u\|^2}{N'(J)};
\end{equation}
where $\widetilde\mu(u)|_{\g{sp}(V_\0)}$ is the restriction of
$\widetilde\mu(v_\0)$ to $\g{sp}(V_\0)$.

\medskip

We put for $X''\in  \big(\g{so}(V_\1)\oplus\g{spo}(V)_\1\big)_\PC$:
\begin{equation}
\phi(X'',v)=\int_{\g{sp}(V_\0)}\r(X')\exp\big(-\frac\ii 2 B(v,(X'+X'')v)\big)
\end{equation}

\begin{lemme} For any $X''\in \big(\g{so}(V_\1)\oplus\g{spo}(V)_\1\big)_\PC$, 
$\phi(X'',v)$ is a well defined rapidly decreasing function on $V$. Moreover,
$\phi$ is polynomial in $X''$.
\end{lemme}

\begin{proof}

Let  $v_\0\in V_\0\tens V_\0^*$ be the generic point of $V_\0$ and $v_\1\in
V_\1\tens V_\1^*$ be the generic point of $V_\1$. As in section \ref{Spf_P} we
have $v=v_\0+v_\1$. Then:
\begin{equation}
\begin{split} B(v,(X'+X'')v)&=B(v_\0,X'v_\0)+B(v,X''v)\\
&=-2\widetilde\mu(v_\0)(X')+B(v,X''v).
\end{split}
\end{equation}

It follows
\begin{equation}
\begin{split}
\phi(X'',v)&=\int_{\g{sp}(V_\0)}\r(X')\exp\big(\ii  \widetilde\mu(v_\0)(X')\big)
\exp\big(-\frac\ii 2 B(v,X''v)\big)\\
&=\what\r\big(-\widetilde\mu(v_\0)|_{\g{sp}(V_\0)}\big)\exp\big(-\frac\ii2
B(v,X''v)\big).
\end{split}
\end{equation}

\medskip

Since $\bb(X'')\in\g{so}(V_\1)$, we have:
\begin{equation}
B(v,X''v)=B(v_\1,\bb(X'')v_\1)+B(v,(X''-\bb(X''))v).
\end{equation}
Hence, $B(v,X''v)$ is nilpotent and $X''\mapsto\exp(-\frac\ii2 B(v,X''v))$
defines a polynomial function on $\g{so}(V_\1)\oplus \g{spo}(V)_\1$ with values
in $S(V^*)$. In particular $\phi$ is polynomial in $X''$.

\medskip

Since $\what \r$ is rapidly decreasing on $\g{sp}(V_\0)^*$, formula
(\ref{eq:MinN}) ensures that $\what\r(-\widetilde\mu(v_\0)|_{\g{sp}(V_{\0})})$
is a rapidly decreasing function on $V_\0$.

Finally, for any $X''\in \big(\g{so}(V_\1)\oplus\g{spo}(V)_\1\big)_\PC$,
$\phi(X'',v)$ is a rapidly decreasing function on $V$. 
\end{proof}

Now, the integral:
\begin{equation}
h(X'')=\int_Vd_V(v)\phi(X'',v)=\int_Vd_V(v)\int_{\g{sp}(V_\0)}\r(X')\exp\big%
(-\frac\ii 2 B(v,(X'+X'')v)\big)
\end{equation}
converges and it defines a polynomial function on
$\g{so}(V_\1)\oplus\g{spo}(V)_\1$. This means that $\Spf$ is a generalized
function on $\g{sp}(V_\0)$ with values in
$S\big((\g{so}(V_\1)\oplus\g{spo}(V)_\1)^*\big)$.

\medskip

\subsubsection{Example: Symplectic $2$-dimensional vector space:}
\label{Spf:gene:Ex} This is the crucial example.
 The problem is to show that for a smooth compactly supported distribution $\r$
on $\g{sp}(V)$:
\begin{equation}
v\mapsto \int_{\g{sp}(V)}\r(X)\exp(-\frac\ii 2 B(v,Xv))
\end{equation}
is a rapidly decreasing function on $V$.

We use notations and results of  sections \ref{2symp} and \ref{2symp2}. We
denote by $(\boldsymbol{a},\boldsymbol{b},\boldsymbol{c})$ the basis of
$\g{sp}(V)^*$ such that $\boldsymbol{a}(X)=a,\,
\boldsymbol{b}(X)=b$ and $\boldsymbol{c}(X)=c$.

Thus:
\begin{equation}
\label{eq:FourierSL2(1)}
\int_{\g{sp}(V)}\r(X)\exp(\ii\mu(X,v))=\what\r\Big(\frac12
\big(\boldsymbol{c}x^2-2\boldsymbol{a}xy-\boldsymbol{b}y^2\big)\Big).
\end{equation}
Since $\what\r$ is a rapidly decreasing function on $\g{sp}(V)^*$,
 the above function is rapidly decreasing on $V$.

\subsection{Comparison with the analytic version of section \ref{analytic}} In
this section we denote by $\Spf_{an}$ the analytic superPfaffian defined on
$\VC^+\times\ii\,\g{spo}(V)$ by formula (\ref{eq:Spf}) and by $\Spf_{gene}$ the
generalized superPfaffian defined on $\g{spo}(V)$ by formula (\ref{eq:SpfII}).

\medskip

Since $X+\ii Y\in\g{spo}(V)\times\ii \VC^-$ is equivalent $\ii(X+\ii
Y)\in\VC^+\times\ii\g{spo}(V)$,  \begin{equation}
(X,\ii Y)\mapsto\ii^{\frac{m-n}2}\Spf_{an}(\ii(X+\ii Y))
\end{equation}
 is an analytic function on $\g{spo}(V)\times\ii\VC^-$. We consider it as an
analytic function on the open cone $
\g{spo}(V)_\0\times\ii\VC^-$ of $\g{spo}(V\tens \mathbb{C})_\0$ with values in
$\Lambda(\g{spo}(V)_\1^*)$. 

We fix a relatively compact open neighborhood $\XC$ of $0$ in $\g{spo}(V)_{\0}$. Since $\Spf$ is homogeneous of degree $\frac{n-m}2$ it follows that for any
relatively compact open subset $\WC\subset\g{spo}(V)_\0$, there exists a
constant $K_\WC$ such that  for any  $(X,\ii Y)\in \WC\times\ii(\VC^-\cap\XC)$ and any
homogeneous differential operator $\DC\in\Lambda(\g{spo}(V)_\1)$ we have for
some $k\in\mathbb{N}$:
\begin{equation}
\Bigl|\bigl(\DC\,\Spf_{an}\big)(\ii(X+\ii
Y))\Big|=\Bigl|\bigl(\DC\,\Spf_{an}\big)(-Y+\ii X))\Big|\leqslant K_\WC
N'(Y)^{-k}.
\end{equation}
($N'$ is a norm on $\g{spo}(V)_\0$.) 

Then \cite[Theorem 3.1.15]{Hor83}   shows that its limit when $Y$ goes to $0$ in
$\VC^-$  exists as a generalized function on $\g{spo}(V)_\0$ with values in
$\Lambda(\g{spo}(V)_\1^*)$. We have:
\begin{equation}
\label{eq:limit}
\Spf_{gene}(X)=\lim_{Y\rightarrow 0,\\
Y\in\VC^-}\ii^{\frac{m-n}2}\Spf_{an}(\ii(X+\ii Y)).
\end{equation}

\medskip

Since $\Spf_{an}$ is the holomorphic extension of $\Spf_{an}|_{\VC^+}$. It is
entirely determined by $\Spf_{an}|_{\VC^+}$. On the other hand, since
$\Spf_{gene}(X)$ is the limit of $\ii^{\frac{n-m}{2}}\Spf_{an}(\ii(X+\ii Y))$,
$\Spf_{gene}$ is determined by $\Spf_{an}$ and thus by $\Spf_{an}|_{\VC^+}$.

\medskip In particular it follows that $\Spf_{gene}$ possesses  the properties
of the sections \ref{Inv}-\ref{homog}.

\medskip

Let $\PC$ be a near superalgebra and $X\in\g{spo}(V)_\PC(\VC^+)$, let $\e>0$, we
have:
\begin{equation}
\label{eq:e}
\Spf_{an}((\e+\ii)X)=(\ii+\e)^{\frac{n-m}2}\Spf_{an}(X)
\end{equation}
It follows, taking the limit of (\ref{eq:e})  when $\e$ goes to zero multiplied
by $\ii^{\frac{m-n}2}$:
\begin{equation}
\Spf_{gene}|_{\VC^+}=\Spf_{an}|_{\VC^+}.
\end{equation}
($\phi|_{\VC^+}$ denotes the restriction of the (generalized) function $\phi$ to
the open set $\VC^+$.) In particular, $\Spf_{gene}$ is analytic on $\VC^+$.

From now on $\Spf$ stands for $\Spf_{gene}$, and for
$(X,Y)\in\g{spo}(V\tens\mathbb{C})_\PC\big(\VC^+\times\ii
\g{spo}(V)\big)$ $\Spf(X+\ii Y)$ stands for $\Spf_{an}(X+\ii Y)$.

\subsection{Evaluation of $\Spf$ on $\VC_{p,q}$}\label{EvalVpq} We recall that
$\UC_{p,q}\subset\g{sp}(V_\0)$ denote the open set of $X\in\g{sp}(V_0)$ such
that $v\mapsto B(v,Xv)$ is a quadratic form of signature $(p,q)$ on $V_{\0}$ and
$\VC_{p,q}=\UC_{p,q}\times\g{so}(V_\1)$.

\begin{prop}\label{prop:EvalVpq}
 Let $\PC$ be a near superalgebra. Let $(p,q)\in\mathbb{N}^2$ such that $p+q=m$.
 Let $X\in\g{spo}(V)_\PC(\VC_{p,q})$. It means that $X$ is represented in a
symplectic basis $(e_1,\dots,e_m,f_1,\dots,f_n)$ by
 \begin{equation}
X=
\begin{pmatrix}A&B\\
C&D
\end{pmatrix}\in \g{spo}(V)_\PC,\ \text{with }\bb(A)\in\UC_{p,q}
\end{equation}
We have:
\begin{equation}
\label{eq:SpfVpq}
\Spf
\begin{pmatrix}A&B\\
C&D
\end{pmatrix} =\ii^{q}\frac{\Pfaff(D-CA^{-1}B)}{\sqrt{\big|\det(A)}\big|}.
\end{equation}
where we recall that $\Pfaff$ is the ordinary Pfaffian.
\end{prop}

\begin{proof} The first equality in formula (\ref{eq:SpfDecomp}) implies that it is enough to
 prove
the formula for $X\in\big(\g{spo} (\VC_{p,q})_{\0}\big)_\PC$.

\medskip

First, consider the particular case where  $X\in\g{spo}(V)_\PC(\VC^+)$. In this
case $(p,q)=(m,0)$, thus the coefficient is $\ii^0=1$ and for
$X\in\g{spo}(V)_\PC(\VC^+)$, $\bb(A)\in\UC^+$ and thus $\det(\bb(A))>0$. The
proposition reduces in this case to formula  (\ref{eq:SpfDecomp}).
   
\medskip

Since $\VC_{p,q}$ is a purely even real manifold, it is enough to consider 
$\PC=\mathbb{R}$ and $X\in\VC_{p,q}$. We put $X=
\begin{pmatrix} A      &   0 \\
   0   &  D
\end{pmatrix},$ with $A\in\UC_{p,q}$ and $D\in\g{so}(V_\0)$.

Then:
\begin{equation}
\Spf(X)=\ii^{\frac{m-n}{2}}\int_{V_\0}d_{V_\0}(v_\0)\exp(-\frac\ii2
B(v_\0,Av_\0))\,
\int_{V_\1}d_{V_\1}(v_\1)\exp(-\frac\ii2 B(v_\1,Dv_\1)).
\end{equation}

On one hand:
\begin{equation}
\int_{V_\1}d_{V_\1}(v_\1)\exp(-\frac\ii2 B(v_\1,Dv_\1))=\ii^{\frac n2}\Pfaff(D).
\end{equation}

On the other hand, it is well known (cf. for example \cite[formula
3.4.6]{Hor83}) that for $A\in\UC_{p,q}$:
\begin{equation}
\int_{V_\0}|dx^1\dots dx^m|\exp\big(-\frac{\ii}{2}B(v_\0,Av_\0)\big)=
\frac{(2\pi)^{\frac m2}}{\exp(\ii\frac {p-q}4\pi)\sqrt{\big|\det(A)\big|}}.
\end{equation}

Since $p+q=m$, $\frac{\ii^{\frac m2}}{\exp(\ii\frac {p-q}4\pi)}=\ii^{ q}$ and
the formula follows.
\end{proof}

\medskip

In particular, it implies that $\Spf$ is smooth on $\VC$ (in fact it is
analytic) and  that for any 
$X\in\g{spo}(V)_\PC(\VC)$, we have:
\begin{equation}
\Spf(X)^2=\Ber^-(X).
\end{equation}

\subsection{Example: $\g{spo}(2,2)$}

Let us consider as an example the case $\g g=\g{spo}(2,2)$. Let $\PC$ be a near
superalgebra. The algebra $\g{spo}(2,2)_\PC$ is the set of matrices:

\begin{equation}
X=
\begin{pmatrix}
\vbox{\offinterlineskip\halign{\hglue 5pt#\hglue 5pt&#&
\hglue 5pt#\hglue 5pt\cr $A$ &\vrule height 12pt depth 10pt& $B$\cr
\noalign{\hrule} $C$&
\vrule height 12pt depth 8pt&$D$\cr}}
\end{pmatrix} =
\begin{pmatrix}
\vbox{\offinterlineskip\halign{\hglue 5pt#\hglue 5pt&#&
\hglue 5pt#\hglue 5pt\cr $
\begin{matrix}a&\hfill b\\
c&-a
\end{matrix}$ &\vrule height 18pt depth 12pt& $
\begin{matrix}\be&\d\cr -\a&-\gamma
\end{matrix}$\cr
\noalign{\hrule} $
\begin{matrix}-\a&\hskip 6pt-\be\cr -\gamma&\hfill-\d\cr
\end{matrix}$\hglue 6pt&
\vrule height 18pt depth 13pt&\hglue 11pt$
\begin{matrix}0&-d\cr  d&0\cr
\end{matrix}$\cr}}
\end{pmatrix}
\end{equation} 
 
 where $a,b,c,d\in\PC_\0$ and $\a,\be,\gamma,\d\in\PC_\1.$ Here $V=\mathbb
R^{(2,2)}$ is endowed with the symplectic form $\bB$ given in the canonical base 
$(e_1,e_2,f_1,f_2)$ ($|e_i|=0$ and $|f_i|=1$) by $\bB(f_i,f_j)=\d_i^j,$
$\bB(e_1,e_2)=-\bB(e_2,e_1)=1,$ $\bB(e_1,e_1)=\bB(e_2,e_2)=0$ and
$\bB(e_i,f_j)=\bB(f_j,e_i)=0$, and  the orientation defined by the basis
$(e_{1},e_{2})$.

\medskip

We have:
\begin{equation}
\label{eq:product}
\begin{split} CA^{-1}B &=\frac{-1}{a^2+bc}
\begin{pmatrix}
\a      &   \be \\
   \gamma   &\d  
\end{pmatrix}
\begin{pmatrix}
      a&   b \\
      c&-a  
\end{pmatrix}
\begin{pmatrix}
\be &\d    \\
-\a      &  -\gamma
\end{pmatrix}\\
&=\frac{1}{a^2+bc}
\begin{pmatrix} 0      &  -(  a\a\d-b\a\gamma+c\be\d+a\be\gamma)\\
 a\a\d-b\a\gamma+c\be\d+a\be\gamma   &  0
\end{pmatrix}.
\end{split}
 \end{equation}

We denote by $\Spf_{0}$ the superpfaffian on $\g{sp}(\mathbb
R^2)$. Thus, for $X\in\g{spo}(2,2)_\PC(\VC)$ ($\VC=\UC\times\g{so}(2)$), we get
from (\ref{eq:SpfVpq}):
\begin{equation}
\label{eq:Spf22}
\Spf (X)=\Big(  d-\frac{a(\a\d+\be\gamma)-b\a\gamma+c\be\d}{a^2+bc}  \Big)
\Spf_{0}
\begin{pmatrix}
     a &b    \\
      c&  -a
\end{pmatrix}.
\end{equation}

\medskip

With (cf. formula (\ref{eq:SpfVpq})):
\begin{eqnarray}
\Spf_{0} 
\begin{pmatrix}
     a &b    \\
      c&  -a
\end{pmatrix}& = &  \frac{1}{\sqrt{-(a^2+bc)}}  \qquad\text{ if }
\begin{pmatrix}
     a & b   \\
      c&  -a
\end{pmatrix}\in \UC_{2,0}\\
\Spf_{0} 
\begin{pmatrix}
     a &b    \\
      c&  -a
\end{pmatrix}& = &  \frac{-1}{\sqrt{-(a^2+bc)}}  \qquad\text{ if }
\begin{pmatrix}
     a &b    \\
      c&  -a
\end{pmatrix}\in U_{0,2}\\
\Spf_{0} 
\begin{pmatrix}
     a &b    \\
      c&  -a
\end{pmatrix}& = &  \frac{\ii}{\sqrt{a^2+bc}}  \qquad\text{ if }
\begin{pmatrix}
     a &b    \\
      c&  -a
\end{pmatrix}\in \UC_{1,1}.
\end{eqnarray}

\medskip

Now we give formula (\ref{eq:restrg}) in this particular case. We denote by
$(x,y,\xi,\eta)$ the system of coordinates on $\mathbb R^{(2,2)}$ dual of 
$(e_1,e_2,f_1,f_2)$.  Let $v=e_1x+e_2y+f_1\xi+f_2\eta$ be the generic point of
$\mathbb{R}^{(2,2)}$. 
We have:

\begin{equation}
\mu(X,v)=d\eta\xi+\a\xi x+\be\xi y+\gamma\eta x+\d\eta y-{\frac
c2}x^2+axy+{\frac b2}y^2.
\end{equation}

Thus, since $d_V(v)=\frac1{2\pi}|dx\,dy|\frac{\de}{\de\xi}\frac{\de}{\de\eta}$:

\begin{multline}
\int_{\mathbb R^{(2,2)}}d_V(v)\exp(\ii\mu(X,v))=\\
\frac{1}{2\pi}\int_{\mathbb R^2}|dx\,dy|\,\Big(\ii d-(\a x+\be y)(\gamma x+\d y)
\bigr)\Bigr)
 \exp \ii(-{\frac c2}x^2+axy+{\frac b2}y^2),
\end{multline}

\medskip

Now we compute the integral on the right hand side. We put 
\begin{equation}
\HC =
\begin{pmatrix}
\frac{\a\d+\be\gamma}2&\be\d   \\
 -  \a\gamma   &  -\frac{\a\d+\be\gamma}2
\end{pmatrix}.
\end{equation}
With this notation we have if  $v_\0=e_1x+e_2y$ is the generic point of
$\mathbb{R}^2$:
\begin{equation}
-(\a x+\be y)(\gamma x+\d y)=\bB( v_\0,\HC v_\0)
\end{equation}

Thus
\begin{equation}
\begin{split} -\int_{\mathbb R^2}|dx\,dy|\,(\a x+\be y)(\gamma x+\d y) 
 &\exp \ii(-{\frac c2}x^2+axy+{\frac b2}y^2)\\
&=\int_{\mathbb R^2}|dx\,dy|\,B( v_\0,\HC v_\0)
 \exp-\frac \ii2    \bB(v_\0,Av_\0)\\
&=2\ii\de_{\HC}\int_{\mathbb R^2}|dx\,dy|\,
 \exp-\frac \ii2    \bB(v_\0,Av_\0).
 \end{split}
\end{equation}

Finally, since
\begin{equation}
\Spf_{0}
\begin{pmatrix}
     a&b   \\
      c&  -a
\end{pmatrix} =\frac{\ii}{2\pi}
\int_{\mathbb{R}^2}|dx\,dy| \exp-\frac \ii2    B(v_\0,Av_\0),
\end{equation}
 we have in $\CC^{-\oo}_{\g{spo}(V)}(\g{spo}(V)_{\0})$:
\begin{equation}
\Spf(X)=\Big(  d+2\Big(\frac{\a\d+\be\gamma}2\frac{\de}{\de
a}+\be\d\frac{\de}{\de b}-\a\gamma\frac{\de}{\de c}\Big) 
\Big)
\Spf_{0}
\begin{pmatrix}
     a &b    \\
      c&  -a
\end{pmatrix}.
\end{equation}

From this equation we deduce again (\ref{eq:Spf22}).
\medskip

\subsection{Singularities and wave front set}\label{WF}

Let  $\singsupp(\Spf)$ be the set of singularities of the superPfaffian.

\begin{lemme} We have
\begin{equation}
\singsupp(\Spf)=\g {spo}(V)_\0\setminus \VC.
\end{equation}
\end{lemme}

\begin{proof}[Proof] Proposition \ref{prop:EvalVpq} implies that  $\Spf$ is 
analytic on $\VC$. Moreover, if $X\in\g{spo}(V)_\0\setminus\VC$ and $\WC$ is a
neighborhood of $X$, $\Spf$ is not bounded on $\WC\cap\VC$. Hence, $X$ is a
singularity.
\end{proof}

\medskip

Let $\phi$ be a generalized function function on $\g{spo}(V)$. We denote by
$\Sigma(\phi)$ the cone in $\g{spo}(V)^*$ defined by the following  (cf.
\cite[Formula (8.1.1)]{Hor83}).

Let $N$ be the norm on $\g{spo}(V)$ defined in section \ref{Spf:gene}. Let 
$f\in\g{spo}(V)^*$ then  (cf. \cite[Formula (8.1.1)]{Hor83})   $f\not\in
\Sigma(\phi)$ if and only if $f=0$ or if  for any $k\in\mathbb{N}$, there exists
$C_k\in\mathbb{R}$ such that for any $h$ in some conic neighborhood of $f$:
\begin{equation}
\big|\what{\phi}(h)\big|\leqslant \frac{C_k}{(1+N(h))^k}
\end{equation}

Let $X\in\singsupp(\Spf)$. We denote by $WF_X(\Spf)$ the wave front set of
$\Spf$ at $X$:
\begin{equation}
\label{def:WF_X} WF_X(\Spf)=\bigcap_{\phi,\,\phi(X)\not=0}\Sigma(\phi \Spf),
\end{equation}
where $\phi$ run in the set of compactly supported smooth distributions on
$\g{spo}(V)_\0$.

We recall that $\widetilde\mu:V\mapsto\g{spo}(V)^*$, and
$\widetilde\mu(v)(X)=\mu(X,v)$, thus
$\widetilde\mu(V_\0)\subset\g{spo}(V)_\0^*$.

\begin{prop}
\begin{equation}
\label{eq:WF_X} WF_X(\Spf)\subset\widetilde\mu(V_\0)\setminus\{0\}.
\end{equation}

\end{prop}

\begin{proof}

It follows from formula (\ref{eq:SpfDecomp}) that the wave front set of $\Spf$ on
$\g{spo}(V)$ at $X$ is equal to the wave front set of $\Spf_{0}$ on $\g{sp}(V_\0)$ at $X$. (The
later is a subset of $\g{sp}(V_\0)\oplus\g{sp}(V_\0)^*$ which is canonically
embedded in $\g{spo}(V)_\0\oplus\g{spo}(V)_\0^*$ by mean of the decomposition
$\g{spo}(V)_\0=\g{sp}(V_\0)\oplus\g{so}(V_\1)$. Thus from now on we assume that
$V=V_\0$.

\medskip

Let $\phi$ be a smooth compactly supported distribution on $\g{sp}(V)$. Then,
$\phi\Spf$ is a compactly supported distribution on $\g{sp}(V)$.

\medskip Let us  precise $\what{\phi\Spf}(h)$.
\begin{equation}
\begin{split}
\what{\phi\Spf}(h)&=\int_{\g{sp}(V)}\phi(X)\Spf(X)\exp\big(-\ii h(X)\big)\\
&= \ii^{\frac
m2}\int_Vd_V(v)\int_{\g{sp}(V)}\phi(X)\exp\big(-\frac{\ii}{2}B(v,Xv)\big)
\exp\big(-\ii h(X)\big)\\
&= \ii^{\frac m2}\int_Vd_V(v)\int_{\g{sp}(V)}
\phi(X)\exp\big(-\ii(h-\widetilde\mu(v))(X))\big)\\
&=\ii^{\frac m2} \int_Vd_V(v)\what \phi\big(h-\widetilde\mu(v)\big)
\end{split}
\end{equation}

Since $\phi$ is smooth and compactly supported, $\what\phi$ is rapidly
decreasing. Thus for any $k\in \mathbb{N}$, there is a constant $K_k$ such that:
\begin{equation}
\big|\what \phi\big(h-\widetilde\mu(v)\big)\big|\leqslant
\frac{K_k}{\big(1+N(h-\widetilde\mu(v))\big)^k}
\end{equation}

 For $h\in\g{sp}(V)^*$ we put
\begin{equation}
d(h)=\min_{v\in V}\{N(h-\widetilde\mu(v))\}.
\end{equation}
 We have if $h\not=0$:
\begin{equation}
d(h)=N(h)d\Big(\frac{h}{N(h)}\Big)
\end{equation}
Moreover, since $\widetilde\mu (V) $ is closed there is $v_h\in V$ such that 
$N\big(h-\widetilde\mu(v_h)\big)=d(h)$ (this propriety determines $v_h$ up to
multiplication by $\pm1$).  Then we have: $
N\big(h-\widetilde\mu(v)\big)\geqslant N\big(h-\widetilde\mu(v_h)\big)=d(h) $
and thus 
$
 N\big(\widetilde\mu(v_h)-\widetilde\mu(v)\big)\leqslant
N\big(h-\widetilde\mu(v_h)\big) +N\big(h-\widetilde\mu(v)\big)\leqslant
2N\big(h-\widetilde\mu(v)\big). $ It follows
\begin{equation}
\begin{split}
\big(1+N\big(h-\widetilde\mu(v)\big)\big)^2 &\geqslant(1+d(h))\big(1+\frac
12N\big(\widetilde\mu(v_h)-\widetilde\mu(v)\big)\big)\\
&\geqslant \frac
12(1+d(h))\big(1+N\big(\widetilde\mu(v_h)-\widetilde\mu(v)\big)\big)
\end{split}
\end{equation}
Therefore:
\begin{equation}
\big|\what \phi\big(h-\widetilde\mu(v)\big)\big|\leqslant \frac{2^{\frac
k2}K_k}{(1+d(h))^{\frac k 2}
\big(1+N\big(\widetilde\mu(v_h)-\widetilde\mu(v)\big)\big)^{\frac k2}}
\end{equation}
Thus:
\begin{equation}
\Big|\what{\phi\Spf}(h)\Big|\leqslant
\frac{2^{\frac k2}K_k}{(1+d(h))^{\frac k2}}
\int_V\frac{d_V(v)}{\big(1+N\big(\widetilde\mu(v_h)-\widetilde\mu(v)\big)\big%
)^{\frac k2}}
\end{equation}

We have:
\begin{equation}
\begin{split} 1+N(\widetilde\mu(v_h)-\widetilde\mu(v))&\geqslant
{1+\big|N(\widetilde\mu(v))-N(\widetilde\mu(v_h))\big|}\\
&\geqslant
\begin{cases}
 1+\frac12   N(\widetilde\mu(v))  & \text{ if  }N(\widetilde\mu(v))\geqslant
2N(\widetilde\mu(v_h)) , \\
1      & \text{otherwise}.
\end{cases}
\end{split}
\end{equation}
\medskip

We recall from section \ref{Spf:gene} that we fixed some $J\in\UC^+$ and put for
$v\in V$ $\|v\|^2=\frac12 B(v,Jv)=-\mu(J,v)$. Then, from formula (\ref{eq:MinN})
we obtain $\|v\|\leqslant\sqrt{N(\widetilde\mu(v))N'(J)}$.

We recall that for $r>0$,  with $m=\dim(V)$:
\begin{equation}
\int_{\big\{v\in V,\|v\|\leqslant r\big\}}d_V(v) =\frac1{(2\pi)^{\frac
m2}}\frac{\pi^{\frac m2}r^m}{\frac m2!}
\end{equation}

We have for $k\geqslant m\geqslant 2$:

\begin{equation}
\begin{split}
\int_V
\frac{d_V(v)}{\big(1+N\big(\widetilde\mu(v_h)-\widetilde\mu(v)\big)\big)^{\frac
k2}} &\leqslant  \int_{\Big\{v\in V,\|v\|\leqslant
\sqrt{2N(\widetilde\mu(v_h))N'(J)}\Big\}}d_V(v)\\
&\hskip 1cm +\int_{V}d_V(v)\frac1{( 1+ \frac12  N(\widetilde\mu(v)) ^{\frac k2}
}\\
&\leqslant 
\frac{\big(2N(\widetilde\mu(v_h))N'(J)\big)^{\frac m2}}{2^{\frac m2}\frac m2!}
+\int_{V}d_V(v)\frac1{( 1+  \frac12 
N(\widetilde\mu(v))^{\frac k2}  }\\
&\leqslant 
\frac{\big(N(\widetilde\mu(v_h))N'(J)\big)^{\frac m2}}{\frac m2!}
+\int_{V}d_V(v)\frac1{( 1+  \frac 1{2N'(J)}\|v\|^2)}\\
\end{split}
\end{equation}
We put $M'_k=2^{\frac k2} K_k\frac{N'(J)^{\frac m2}}{\frac m2!} $ and
$M_{k}''=2^{\frac k2} K_k\int_{V}d_V(v)\frac1{( 1+
 \frac 1{2N'(J)}\|v\|^2)}$. Thus:

\begin{equation}
\Big|\what{\phi\Spf}(h)\Big|\leqslant
\frac{M'_kN(\widetilde\mu(v_h))^{\frac m2}+M_{k}''}{(1+d(h))^{\frac k2}}.
\end{equation}

\medskip

Let $\CC$ be  a conic neighborhood of $f$. We put:
\begin{equation}
M_\CC=\min\big\{d(h)\,/\, h\in\CC\text{ and }N(h)=1\big\}.
\end{equation}
Now assume that $f\not\in\widetilde\mu(V)\setminus\{0\}$ and $f\not=0$
($0\not\in WF_X(\Spf)$ by definition). Since $\widetilde\mu(V)$ is closed, there
is a conic neighborhood $\CC$ of $f$ whose closure do not intersect
$\widetilde\mu(V)\setminus\{0\}$. In this case, $M_\CC>0$.  It follows that for
any $h\in\CC$, $h\not=0$:

\begin{equation}
N(\widetilde\mu(v_{h}))\leqslant d(h)+N(h)=d(h)(1+\frac1{d(\frac h{N(h)})})
\leqslant (1+d(h))(1+\frac1{M_{\CC}})
\end{equation}
Thus
\begin{equation}
\begin{split}
\Big|\what{\phi\Spf}(h)\Big|&\leqslant
\frac{M'_k (1+d(h))^{\frac m2}(1+\frac1{M_{\CC}})^{\frac
m2}+M_{k}''}{(1+d(h))^{\frac k2}}\\
&\leqslant
\frac{M'_k(1+\frac1{M_{\CC}})^{\frac m2}+M_{k}''}{(1+d(h))^{\frac {k-m}2}}.
\end{split}
\end{equation}
We put $M_{k}=M'_k(1+\frac1{M_{\CC}})^{\frac m2}+M_{k}''$. It follows for any
$h\in\CC$:
\begin{equation}
\Big|\what{\phi\Spf}(h)\Big|\leqslant
\frac{M_k}{(1+M_\CC N(h))^{\frac {k-m}2}}.
\end{equation}
This proves that $f\not\in WF_X(\Spf)$.

\end{proof}

\subsection{Uniqueness results}\label{UniquenessResults} We put:
\begin{equation}
(\VC^-)^0=\big\{f\in\g{spo}(V)_\0^*\,/\, \forall X\in \VC^-,f(X)\geqslant
0\big\}.
\end{equation}
We have $\widetilde\mu(V_\0)\subset(\VC^-)^0$.

As a direct application of \cite[Theorem 8.4.15]{Hor83} we obtain.

\begin{theo} Let $V$ be a symplectic supervector space. 

Let $\VC^+$, $\VC^-$,
$(\VC^-)^0$, $\singsupp(\Spf)$, be defined as above.

Let $\phi\in\CC^{-\oo}_{\g{spo}(V)}(\g{spo}(V)_{\0})$ be a generalized function
on  $\g{spo}(V)$, such that:
\begin{enumerate}
  \item\label{i} $\phi$ is smooth on $\VC^+$ and $(\phi|_{\VC^+})^2=Ber^-$
  \item \label{ii}$WF(\phi)\subset \g{spo}(V)_\0\times (\VC^-)^0$.
\end{enumerate}

Then, $\phi=\Spf$ or $-\Spf$. More precisely, an orientation of $V_\1$ choose
between $\Spf$ and $-\Spf$.

\end{theo}

\begin{proof}

Let $\phi\in\CC^{-\oo}_{\g{spo}(V)}(\g{spo}(V)_{\0})$ satisfying the above
conditions. Condition (\ref{ii}) and
\cite[Theorem 8.4.15]{Hor83} imply that for any open convex cone $\C$ with
closure included in $\VC^-\bigcup\{0\}$, there is an analytic function $F$ on
$\g{spo}(V)\times\ii\C$ such that on $\g {spo}(V)$:
\begin{equation}
\phi(X)=\lim_{Y\rightarrow 0,\,Y\in\C}F(X+\ii Y).
\end{equation}

We choose an orientation of $V_{\1}$. Condition (\ref{i}) implies that  on
$\VC^+$, $\phi(X)=\pm\Spf(X)$. Thus, by
\cite[Theorem 3.1.15 and Remark]{Hor83}, we have on $\g{spo}(V)\times\ii\C$ 
\begin{equation}
\label{eq:analext} F(X+\ii Y)=\pm\ii^{\frac{m-n}2}\Spf(\ii(X+\ii Y)).
\end{equation}
where we recall from section \ref{analytic} that $\Spf(\ii(X+\ii Y))$ is the
analytic function on $\VC^+\times\ii\,\g{spo}(V)$ defined by formula
(\ref{eq:Spf}). Thus equation  (\ref{eq:limit}) implies that $\phi=\pm\Spf$.

\medskip

Since we saw that  changing the orientation of $V_\1$ changes $\Spf$ to $-\Spf$,
the last remark follows.

\end{proof}

\emph{Remark:} In the definition of $\Spf$ besides the orientation of $V_\1$, we
chose a square root $\ii$ of $-1$. Changing $\ii$ into $-\ii$ changes $\Spf(X)$
into $\overline\Spf(X)$ which is the limit of
$(-\ii)^{\frac{m-n}{2}}\Spf(-\ii(X+\ii Y))$ when $Y$ goes to $0$ in
$\VC^+=-\VC^-$. In particular:
\begin{equation}
WF(\overline\Spf)=-WF(\Spf).
\end{equation}

\medskip

\subsection{Linear subspaces and Subalgebras}\label{Subalgebras}

Let $\g g$ be a linear subspace of $\g{spo}(V)$. We put:
\begin{equation}
\g g_\0^\dag=\big\{f\in\g {spo}(V)_{\0}^* \,/\, f(\g g_{\0})=\{0\}\big\}
\end{equation}

\medskip

 Let us assume that $\g g_{0}^\dag\cap\widetilde\mu(V_{\0})=\{0\}$. This means:
\begin{equation}
\forall u\in V_0\setminus\{0\},\,\exists X\in \g
g_{\0},\,\mu(X,u)=-\frac12B(u,Xu)\not= 0.
\end{equation}

In particular this condition is implied by the following:
\begin{equation}
\g g_{\0}\cap\VC^-\not=\emptyset.
\end{equation}

\medskip

Indeed, assume that $\g g_{\0}\cap\VC^-\not=\emptyset$, then there exists $X\in\g g_{\0}$ 
such that $u\mapsto B(u,Xu)$ is negative definite on $V_{\0}$. in particular for any 
$u\in V_{\0}\setminus\{0\}$, $\mu(X,u)>0$ and
 $\widetilde \mu(u)\not\in\g g_{\0}^{\perp}$.

\medskip When this condition is realized, by standard results on generalized
functions (cf. for example \cite[Corollary 8.2.7]{Hor83}), we can define
$\Spf|_{\g g}$. It is defined for a smooth compactly supported distribution $t$
on $\g g$ 
by:
\begin{equation}
\label{eq:restrg}
\int_{\g g}t(X)\Spf(X)=\ii^{\frac {m-n}2}\int_Vd_V(v)\int_{\g
g}t(X)\exp\big(-\frac\ii2B(v,Xv)\big).
\end{equation}
Let us show that the preceding formula is  meaningful.

\medskip

Again, we take the notations of section \ref{Spf:gene}. Since $\g g_\0$ is a
subalgebra of $\g{spo}(V)_\0$ and $\g{sp}(V_\0)$ and $\g{so}(V_\1)$ are simple
we have $\g g_\0=(\g g_\0\cap\g{sp}(V_\0))\oplus(\g g_{\0}\cap\g{so}(V_\1))$.
Here $\r$ will be a smooth compactly supported distribution on $\g
g_\0\cap\g{sp}(V_\0)$ and $X''\in\big((\g g_\0\cap\g{so}(V_\1))\oplus\g
g_\1\big)_\PC$. We have ($v$ (resp. $v_\0$) is the generic point of $V$ (resp.
$V_\0$)):

\begin{equation}
\label{eq:phi}
\phi(X'',v)=\what\r\big(-\widetilde\mu(v_\0)|_{\g g_\0\cap\g{sp}(V_\0)}\big)
\exp\big(-\frac\ii2B(v,X''v)\big).
\end{equation}

We put:
\begin{equation}
p=\min \{ N(\widetilde\mu(u)|_{\g g_\0\cap\g{sp}(V_\0)}) \,/\, u\in V_\0,\ 
{N(\widetilde\mu(u))}=1\}.
\end{equation}
Formula (\ref{eq:MinN}) gives for $v\in V_\0$:
\begin{equation}
N(\widetilde\mu(u)|_{\g g_\0\cap\g{sp}(V_\0)})\geqslant p
N(\widetilde\mu(u))\geqslant \frac{p}{N'(J)}\|u\|^2.
\end{equation}
 The hypothesis implies that for $u\in V_\0\setminus\{0\}$,
$N(\widetilde\mu(u)|_{\g g_\0\cap\g{sp}(V_\0)})>0$. Thus, since
$\big\{\widetilde\mu(u)\,/\, u\in V_\0,\, N(\widetilde\mu(u))=1\}$ is compact,
we have $p>0$. Hence, since $\what\r$ is rapidly decreasing on $\g g^*$,
$\what\r\big(-\widetilde\mu(v_0)|_{\g g_\0\cap\g{sp}(V_\0)}\big)$ is rapidly
decreasing on $V_\0$ and thus formula (\ref{eq:restrg}) is meaningful.

\subsubsection{Symplectic $2$-dimensional vector spaces}

Let $V$ be a Symplectic $2$-dimensional vector space. Let $\g
g\subset\g{spo}(V)\simeq \g{sl}(2)$ be a subalgebra. Then, $\g g=\g g_{\0}$ and
the condition  $\g g_{\0}^\dag\cap\widetilde\mu(V)\setminus\{0\}=\emptyset$
implies that $\g g=\g{sl}(2)$ or $\g g$ is a compact Cartan subalgebra.

\subsection{The case $V=V_\0$: the superPfaffian as Fourier transform of a
coadjoint orbit}\label{Four}

In this section we assume that $V=V_\0$. We recall from formula
(\ref{eq:moment}) that we denote by $\widetilde\mu$ the moment map from $V$ to
$\g{sp}(V)^*$; for any  $v\in V$ and any $X\in\g{sp}(V)^*$:
\begin{equation}
\widetilde\mu(v)(X)=-\frac12 B(v,Xv).
\end{equation}

We denote by $\widetilde\mu(V)$ the image of $V$  by $\widetilde\mu$ in
$\g{sp}(V)^*$. It is the disjoint union of a nilpotent orbit $\W$ and of
$\{0\}$. 

Let $\W=\widetilde\mu(V)\setminus\{0\}$. It is naturally endowed with an
Hamiltonian structure. We denote by $\w_{\W}$ its symplectic form and by
$\mu_\W$ its moment map which is the identity map from $\g{sp}(V)$ onto
$\g{sp}(V)\simeq\big(\g{sp}(V)^*\big)^*\subset\CC^\oo(\W)$.  
The form $\w_{\W}^{\frac m2}$ determines an orientation on $\W$. By definition 
the Fourier transform of $\W$ is  the 
generalized function on $\g{sp}(V)$ defined by:

\begin{equation}
\FC_{\W}(X)
=\frac{1}{(2\pi)^{\frac{m}{2}}}\int_\W\exp\big(\ii(\mu_\W(X)+\w_{\W})\big).
\end{equation}

Let us  consider $\widetilde\mu$ as a morphism of manifolds from
$V\setminus\{0\}$ to $\W.$ Since $B$ is an antisymmetric bilinear form on $V$,
it is an element of $\Lambda(V^*)$ and thus determines a differential form
$\w_B$ of degree $2$ on $V$. The Liouville integral of section \ref{Symp}
satisfies:
\begin{equation}
d_V(v)=\frac{1}{(2\pi)^{\frac m2}\frac{m}2!}\big|\w_B^{\frac{m}{2}}(v)\big|.
\end{equation}
Moreover the form $\w_B^{\frac m2}$ fix an orientation of $V$. Since $\mu$
induces a morphism of Poisson algebras $\check\mu:S(V^*)\rightarrow
\g{sp}(V)\subset\CC^{\oo}(\W)$ (cf. (\ref{eq:Poisson})), we have 
\begin{equation}
\widetilde\mu^*(\mu_\W)=\widetilde\mu
\text{ and }
\widetilde\mu^*(\w_{\W})=\w_B.
\end{equation}

\medskip

Since $\widetilde\mu(V\setminus \{0\})$  is a double cover of $\W$, we have:
\begin{equation}
\begin{split}
Spf(X)&=\ii^{\frac{m}{2}}\int_Vd_V(v)\exp(-\ii\widetilde\mu(X)(v))\\
&=\frac{1}{(2\pi)^{\frac{m}{2}}}\int_{V\setminus\{0\}}
\exp\big(\ii(\widetilde\mu(X)(v)+\w_B(v))\big)\\
&=2\frac{1}{(2\pi)^{\frac{m}{2}}}\int_\W\exp(\ii(\mu_\W(X)+\w_{\W}))\\
&=2\FC_\W(X).
\end{split}
\end{equation}

\subsubsection{Example: Symplectic $2$-dimensional vector space.} 
We take the same notations as in section \ref{Spf:gene:Ex}.

We identify  $X\in\g{sp}(V)\simeq \g{sl}(2)$ with its matrix in the basis
$(e_1,e_2)$. We consider the restriction of $\Spf$ to the compact cartan
subalgebra:
\begin{equation}
\g t=\{
\begin{pmatrix} 0&-c\\
c&0
\end{pmatrix}/c\in\mathbb R\}.
\end{equation}

As a generalized function on $\g t$ we have:
\begin{equation}
\Spf
\begin{pmatrix} 0&-c\\
c&0
\end{pmatrix}=\lim_{\e\to 0,\,\e>0}\frac 1 {c+\ii \e}
\end{equation}

\medskip

On the other hand, Let $\a\in\g t^*$ such that:
\begin{equation}
\a
\begin{pmatrix} 0&-c\\
c&0
\end{pmatrix}=2\ii c.
\end{equation}
Then $(\a,-\a)$ is the root system of $(\g{sl}(2)\tens\mathbb C,\g
t\tens{\mathbb C})$. Let $\W$ be the nilpotent coadjoint orbit of $\g {sl}(2)$
defined by
\begin{displaymath}
\W=\big\{
\begin{pmatrix} a&b\\
c&d
\end{pmatrix}\in\g{sl}(2)\,\big/\, a^2+bc=0,\ c>0\big\}.
\end{displaymath}
 Then
\begin{displaymath}
\widetilde\mu(V)=\W\cup \{0\}.
\end{displaymath}

As it is well known, the restriction of the Fourier transform of $\W$ to regular
elements $H=\begin{pmatrix} 0&-c\\
c&0
\end{pmatrix}\in\g t$ is:
\begin{equation}
\FC_{\W}(H)=\int_{\W}\exp(\ii\mu_{\W}(H)+\w_{\W})=\frac \ii{\a(H)}.
\end{equation}
It follows that:
\begin{equation}
\Spf
\begin{pmatrix} 0&-c\\
c&0
\end{pmatrix}=2\FC_{\W}
\begin{pmatrix} 0&-c\\
c&0
\end{pmatrix}.
\end{equation}


\begin{thebibliography}{BGV92}

\bibitem[Ber87]{Ber87} F.A. Berezin.
\newblock {\em Introduction to Superanalysis}.
\newblock MPAM D. Reidel Publishing Company, 1987.

\bibitem[BGV92]{BGV92} N.~Berline, E.~Getzler, and M.~Vergne.
\newblock {\em Heat Kernels and Dirac Operators}.
\newblock Springer-Verlag, 1992.

\bibitem[Bis86]{Bis86} J.-M. Bismut.
\newblock Localizations formulas, superconnections, and the index theorem for
  families.
\newblock {\em Communications in Mathematical Physics}, 103:127--166, 1986.

\bibitem[BL75]{BL75} F.~A. Berezin and D.~A. Leites.
\newblock Supermanifolds.
\newblock {\em Soviet Math. Dokl.}, Volume 16, n$^0$5:1218--1222, 1975.

\bibitem[BV83]{BV83a} N.~Berline and M.~Vergne.
\newblock Z\'eros d'un champ de vecteurs et classes caract\'eristiques
  \'equivariantes.
\newblock {\em Duke Mathematical Journal}, Volume 50 pages 539--549, 1983.

\bibitem[DeW84]{DeW84} B.~De~Witt.
\newblock {\em Supermanifolds}.
\newblock Cambridge University Press, 1984.

\bibitem[DV88]{DV88} Michel Duflo and Mich\`ele Vergne.
\newblock {\em Orbites coadjointes et cohomologie
\'equivariante}, volume~82 of
  {\em Progress in Mathematics}, pages 11--60.
\newblock Birk{\"a}user, 1988.

\bibitem[GJ81]{GJ81} J.~Glimm, and A.~Jaffe.
\newblock {\em Quantum Physics}.
\newblock Springer-Verlag, 1981.

\bibitem[HC64]{HC64a} Harish-Chandra.
\newblock Invariant distributions on Lie algebras.
\newblock {\em American Journal of Mathematics}, 86:271--309, 1964.

\bibitem[HC65]{HC65a} Harish-Chandra.
\newblock Invariant eigendistributions on a semisimple Lie algebras.
\newblock {\em Publications Math\'ematiques de l'IHES}, 27:5--54, 1965.

\bibitem[H{\"o}r83]{Hor83} L.~H{\"o}rmander.
\newblock {\em The Analysis of Linear Partial Differential Operators I}.
\newblock Springer-Verlag, 1983.

\bibitem[Kos77]{Kos77} B.~Kostant.
\newblock Graded manifolds, graded lie theory and prequantization.
\newblock In {\em LNM 570}, pages 177--306. Springer-Verlag, 1977.

\bibitem[Lav98]{Lav98} P.~Lavaud.
\newblock Formule de localisation en superg\'eom\'etrie.
\newblock Th\`ese de doctorat de l'Universit\'e de Paris VII, 1998.

\bibitem[Lav04]{Lav04} P.~Lavaud.
\newblock Equivariant Cohomology and Localization Formula in Supergeometry
\newblock Preprint

\bibitem[Man88]{Man88} Yu.~I. Manin.
\newblock {\em Gauge Fields and Complex Geometry}.
\newblock Springer-Verlag, 1988.

\bibitem[MQ86]{MQ86}
V.~Matha\"i and D.~Quillen.
\newblock Superconnections, thom classes, and equivariant differential forms.
\newblock {\em Topology}, Volume 25, No 1:85--110, 1986.

\bibitem[Ste90]{Ste90}
J.~Stembridge.
\newblock Non intersecting paths and Pfaffians.
\newblock {\em Adv. in Math.}, Volume 83 96--131, 1990.

\bibitem[Vor91]{Vor91} T.~Voronov.
\newblock Geometric integration theory on supermanifolds.
\newblock In {\em Mathematical Physics Reviews}, Volume 9, Part 1.
  Harwood Academic Publishers, 1991.

\bibitem[Wei53]{Wei53} A.~Weil.
\newblock Th\'eorie des points proches sur les vari\'et\'es diff\'erentiables.
\newblock In {\em Colloque de g\'eom\'etrie diff\'erentielle}, pages 111--117,
  1953.
  

\bibitem[Wey46]{Wey46} H.~Weyl.
\newblock {\em The Classical Groups}.
\newblock Princeton University Press, 1946.

\end{thebibliography}
\end{document}